\documentclass[a4paper,fleqn]{cas-dc}

\usepackage[authoryear,longnamesfirst]{natbib}
\usepackage{mathtools}
\usepackage{amsmath}
\usepackage{amssymb}
\usepackage{amsfonts}
\usepackage{graphicx}
\usepackage{textcomp}
\usepackage{hyperref}
\usepackage{comment}
\usepackage{dsfont}
\usepackage{enumerate}
\usepackage{wrapfig}
\newcommand{\diag}{\mathsf{diag}}
\newcommand{\cN}{{\mathcal{N}}}
\newcommand{\cI}{{\mathcal{I}}}
\newcommand{\cD}{{\mathcal{D}}}
\newcommand{\cG}{{\mathcal{G}}}
\newcommand{\cE}{{\mathcal{E}}}
\DeclareMathAlphabet{\mymathbb}{U}{BOONDOX-ds}{m}{n}
\newcommand{\one}{\mathds{1}}
\newcommand{\zero}{\mymathbb{0}}
\newcommand*{\QEDB}{\hfill\ensuremath{\square}}

\newtheorem{defn}{Definition}
\newtheorem{assump}{Assumption}
\newtheorem{thm}{Theorem}
\newtheorem{rem}{Remark}

\newtheorem{prop}{Proposition}

\newtheorem{ex}{Example}
\newtheorem{com}{Comment}

\def\tsc#1{\csdef{#1}{\textsc{\lowercase{#1}}\xspace}}
\tsc{WGM}
\tsc{QE}

\begin{document}
\let\WriteBookmarks\relax
\def\floatpagepagefraction{1}
\def\textpagefraction{.001}

\shorttitle{Control of Power Grids With Switching Equilibria}    

\shortauthors{Abdelgalil et al.}  

\title [mode = title]{Omega-Limit Sets and Input-to-State Stability in Power Grids With Switching Equilibria}  

\tnotemark[1] 
\tnotemark[2]

\tnotetext[1]{This work was supported in part by the National Science Foundation (NSF) awards 2444163 and 2305756.} 

\tnotetext[2]{This work was authored in part by the National Renewable Energy Laboratory for the U.S. Department of Energy (DOE) under Contract No. DE-AC36-08GO28308. Funding provided by theAdvanced Modeling Grid Research Program at the Office of Electricity, U.S. Department of Energy. The views expressed in the article do not necessarily represent the views of the DOE or the U.S. Government. The U.S. Government retains and the publisher, by accepting the article for publication, acknowledges that the U.S. Government retains a nonexclusive, paid-up, irrevocable, worldwide license to publish or reproduce the published form of this work, or allow others to do so, for U.S. Government purposes.}

%

\author[1]{Mahmoud Abdelgalil}[orcid=0000-0003-1932-5115]
\ead{mabdelgalil@ucsd.edu}
\credit{Conceptualization, Formal analysis, Methodology, Validation, Writing - original draft, Writing - review \& editing}

\author[1]{Vishal Shenoy}[orcid= 0009-0001-7572-1916]

\ead{vshenoy@ucsd.edu}

\cormark[1]

\credit{Conceptualization, Formal analysis, Methodology, Software, Validation, Writing - original draft, Writing - review \& editing}

\author[2]{Guido Cavraro}[orcid= 0000-0003-0296-720X]
\ead{guido.cavraro@nrel.gov}

\credit{Conceptualization, Funding acquisition, Methodology, Writing - review \& editing}

\author[3]{Emiliano Dall'Anese}[orcid= 0000-0002-6486-7477]
\ead{edallane@bu.edu}

\credit{Conceptualization, Funding acquisition, Methodology, Writing - original draft, Writing - review \& editing}

\author[1]{Jorge I. Poveda}[orcid=0000-0003-0733-951X]
\ead{jipoveda@ucsd.edu}

\credit{Conceptualization, Funding acquisition, Methodology, Supervision, Writing - original draft, Writing - review \& editing}

\affiliation[1]{organization={Department of Electrical and Computer Engineering, University of California San Diego},
            country={USA}
            }
\affiliation[2]{organization={Power Systems Engineering Center, National Renewable Energy Laboratory},
            country={USA}
            }
\affiliation[3]{organization={Department of  Electrical and Computer Engineering and Division of Systems Engineering, Boston University},
            country={USA}
            }

\cortext[1]{Corresponding author}



\begin{abstract}
This paper studies a power transmission system with both conventional generators (CGs) and distributed energy assets (DEAs) providing frequency control. We consider an operating condition  with demand aggregating two dynamic components: one that switches between different values on a finite set, and one that varies smoothly over time. Such dynamic operating conditions may result from protection scheme activations, external cyber-attacks, or due to the integration of dynamic loads, such as data centers. Mathematically, the dynamics of the resulting system are captured by a system that switches between a finite number of vector fields ---or modes--, with each mode having a distinct equilibrium point induced by the demand aggregation. To analyze the stability properties of the resulting switching system, we leverage tools from hybrid dynamic inclusions and the concept of $\Omega$-limit sets from sets. Specifically, we characterize a compact set that is semi-globally practically asymptotically stable under the assumption that the switching frequency and load variation rate are sufficiently slow. For arbitrarily fast variations of the load, we use a level-set argument with multiple Lyapunov functions to establish input-to-state stability of a larger set and with respect to the rate of change of the loads. The theoretical results are illustrated via numerical simulations on the IEEE 39-bus test system.
\end{abstract}




\maketitle

The modern power grid integrates a diverse mix of generation resources, ranging from conventional generators (CGs) to distributed energy assets (DEAs). Although DEAs were initially introduced to enhance grid resilience and improve efficiency, they also introduce new challenges. In particular, their performance is affected by changes in prevailing weather conditions, frequency and voltage fluctuations, and increased vulnerability to cyber-attacks \cite{horowitz2019overview, dibaji2019systems}. These concerns have spurred extensive research on the reliability and failure modes of power networks \cite{kavvathas2024resilient, ghiasi2019analytical}.

In the case of failures of frequency-responsive DEAs, each failure mode corresponds to a different net load configuration and different configurations in the overall system inertia and damping. This difference arises from several factors, including the presence of intermittent loads and the implementation of countermeasures against cyber-attacks.

In such scenarios, using the standard swing equation \cite{kundur2007power} to model the dynamics of CGs and frequency-responsive DEAs  results in a linear, time-invariant system that switches between a finite number of stable vector fields, which in general may have different equilibrium points. As shown in \cite{veer2019switched,alpcan2010stability,yin2023stability,baradaran2020omega}, the presence of multiple equilibria precludes the use of standard stability analyses, commonly applied in the frequency control literature, which typically only consider switching between vector fields. This necessitates the development of new tools to study the asymptotic behavior of controllers in power systems with switching dynamically aggregated demand.
 
{\color{black} Stability tools from switching systems' theory have been used to investigate the asymptotic properties of voltage source converters \cite{tang2023common}, leveraging Lyapunov functions to derive estimates of the region of attraction of a stable common equilibrium point. 
In \cite{feng2025online}, the authors design frequency controllers for power networks with time-varying inertia. In this case, each inertia level corresponds to a different mode, resulting in a nonlinear switching system that was studied using standard Lyapunov arguments. In \cite{chen2023distributed}, the authors consider the dual objective problem of simultaneously meeting reference requirements and energy balance in energy storage systems. The theory of switching systems has also been used to study time-varying topologies arising due to unreliable communications between multiple agents in networked systems; voltage regulation in DC microgrids has also been studied under event-triggered switches in \cite{peng2022distributed}.} {\color{black} In \cite{susuki2011hybrid}, the authors use tools from reachability analysis of hybrid automata (a modeling approach equivalent to \eqref{eqn:hy_sys_def} via \cite[Chapter 1, Section 4]{goebel2012hybrid}) to analyse the transient stability of both single and double-machine infinite buses. }

For a recent survey on switching and hybrid systems' tools applied to power systems we refer the reader to \cite{ochoa2023control} and \cite{odunlami2025hybrid}.

In contrast to these works, this paper considers a setting where the controlled power system exhibits switching dynamics without a common equilibrium point. While this scenario is more representative of practical conditions, it has received considerably less attention due to the complexity of the resulting ``stable'' behavior when switching occurs sufficiently slowly. Indeed, to the best of our knowledge, a detailed characterization of this behavior and its stability and robustness properties in power networks with primary controllers and dynamic loads have not been previously addressed in the literature.  

Motivated by the previous background, in this paper we make the following contributions: \emph{(i)} By modeling each operating mode of the power dynamics as an LTI system, we provide an explicit characterization of a compact set that is asymptotically stable (in a semi-global practical sense) for the switching dynamics whenever the switching is sufficiently slow. This compact set is constructed on the basis of the data of the matrices that characterize each of the modes of the system, and by leveraging a modeling framework based on hybrid dynamic inclusions that encodes the modes of the system as logic states, and the switches as events triggered by a dynamic resetting timer. \emph{(ii)} By introducing a novel explicit characterization of $\Omega$-limit sets for switching linear systems, we show that the stability properties of such a set are robust with respect to slow variations of the load aggregation, provided the 2-norm of the difference between any two ``internal'' matrices characterizing the modes is sufficiently small. \emph{(iii)} For arbitrarily fast variations of
the load, we use a class of hybrid Lyapunov functions to establish an input-to-state stability-like result that characterizes the convergence of the trajectories of the system to a larger set that depends on the distance between the different equilibrium points induced by the loads. \emph{(iv)} Finally, we validate the theoretical results via numerical simulations on the IEEE 39-bus test system with both CGs and DEAs. To the best of our knowledge, our results provide the first stability characterization for power transmission systems with switching net loads and DEA failures under dwell time conditions in terms of (semi-globally practically) asymptotically and input-to-state stable sets for switching systems with distinct equilibrium points. 

The rest of this paper is organized as follows. Section \ref{sec:Prelims} presents the preliminaries and notation. In Section \ref{sec:prob_form}, we model the power transmission system as a hybrid system present the main stability results. We present numerical simulations in Section \ref{sec:sims} on the IEEE 39-bus test system. Finally, we conclude in Section \ref{ref:conc}. 

\section{Preliminaries}
\label{sec:Prelims}
 We use $\mathbb{R}^n$ to denote the $n$-dimensional Euclidean space and $\mathbb{Z}_{\geq 0}$ for the set of positive integers. The set of non-negative real numbers is denoted by $\mathbb{R}_{\geq 0}$. We denote the Euclidean norm of $x \in \mathbb{R}^n$ by $\lvert x \rvert$. Given $x \in \mathbb{R}^n$ and a closed set $\mathcal{A} \subset \mathbb{R}^n$, the distance of $x$ to $\mathcal{A}$ is denoted $\lvert x\rvert_{\mathcal{A}}$ and it is defined by $\lvert x \rvert_{\mathcal{A}} \coloneqq \inf_{y \in \mathcal{A}} \lvert x- y \rvert$. Given $r >0$, the set $\{x \in \mathbb{R}^n \, : \, \lvert x \rvert \leq r \}$ is denoted by $r\mathbb{B}$. Given a matrix $A\in\mathbb{R}^{n\times m}$, we use $|A|$ to denote the induced matrix norm of $A$, defined as $|A|:=\sup_{x\in\mathbb{B}}|Ax|$. Given two vectors $u, v$, their concatenation $[u^\top ~ v^\top]^\top$ is denoted as $(u,v)$. The closure of a set $S \subset \mathbb{R}^n$ is denoted by $\Bar{S}$. The closure of the convex hull of $S$ is denoted as $\overline{\text{co}}~S$. For a symmetric matrix $P$, let $\lambda_{\text{min}}(P)$ and $\sigma_{\text{min}}(P)$ denote the minimum eigenvalue and singular value of $P$ respectively. We define $\lambda_{\text{max}}$ and $\sigma_{\text{max}}$ analogously. A set-valued mapping $M:\mathbb{R}^p\rightrightarrows\mathbb{R}^n$ is outer semicontinuous (OSC) at $z$ if for each sequence $\{z_i,s_i\}\to(z,s)\in\mathbb{R}^p\times\mathbb{R}^n$ satisfying $s_i\in M(z_i)$ for all $i\in\mathbb{Z}_{\geq0}$, we have $s\in M(z)$. A mapping $M$ is locally bounded (LB) at $z$ if there exists an open neighborhood $N_z\subset\mathbb{R}^p$ of $z$ such that  $M(N_z)$ is bounded. 

A function $\alpha : \mathbb{R}_{ \geq 0} \to \mathbb{R}_{\geq 0}$ {\color{black} is in class} $\mathcal{K}$ if it is zero at zero, continuous, and strictly increasing.    

A function is {\color{black} in }class $\mathcal{K}_\infty$ if it is in class $\mathcal{K}$ and unbounded.
 
A function $\beta : \mathbb{R}_{\geq 0} \to \mathbb{R}_{\geq 0}$ is {\color{black} in class }$\mathcal{KL}$ if $\beta(\cdot, t)$ is of class $\mathcal{K}$ for each fixed $t \geq 0$ and $t \to \beta(r,t)$ is non-increasing and decreases to zero as $t \to \infty$ for each fixed $r \geq 0$. {\color{black} Let $\mathcal{Q} := \{1,2, \dots, l\}$}. The class of right-continuous, piece-wise constant functions from $\mathbb{R}_{\geq 0} \to \mathcal{Q} \subset \mathbb{Z}_{\geq 0}$ is denoted by $\mathcal{S}$. Also, we denote by $N_{\sigma}(t_1,t_2)$ the number of switches exhibited by a signal $\sigma \in \mathcal{S}$ over the time interval $[t_1,t_2]$. A signal $\sigma \in \mathcal{S}$ {\color{black} satisfies} the average dwell-time (ADT) condition if for any two times { \color{black} in the domain of $\sigma$} we have:
\begin{equation}
    \label{eq:ADT}
    \forall~t_2 \geq t_1 \geq 0,~~~N_{\sigma}(t_1,t_2) \leq N_0 + \frac{t_2-t_1}{\tau_d}.
\end{equation}     
In this paper, we use the formalism of hybrid dynamic inclusions \cite{goebel2012hybrid} to model the systems of interest. A hybrid dynamic inclusion $\mathcal{H}$ is characterized by the tuple $(C,F,D,G)$ and the following expressions:
\begin{subequations}\label{eqn:hy_sys_def}
\begin{align}
&\xi\in C,~~~~~~\dot{\xi} \in F(\xi), \label{flowsHDS}\\
&\xi\in D,~~~~ \xi^+ \in G(\xi),\label{jumpsHDS}
\end{align}
\end{subequations}
where $\xi \in \mathbb{R}^n$ is the state, $C\subset \mathbb{R}^n$ is the \emph{flow set}, $D \subset \mathbb{R}^n$ is the \emph{jump set}, $F : \mathbb{R}^n \rightrightarrows \mathbb{R}^n$ is the \emph{flow map}, and $G : \mathbb{R}^n \rightrightarrows \mathbb{R}^n$ is the \emph{jump map}. Since system \eqref{eqn:hy_sys_def} combines continuous-time and discrete-time dynamics, its solutions $\xi(t,j)$ are parametrized by a continuous-time index $t\in\mathbb{R}_{\geq0}$, which increases continuously during the flows \eqref{flowsHDS}, and a discrete-time index $j$, which increases by one during the jumps \eqref{jumpsHDS}. Therefore, solutions $\xi(t,j)$ to system \eqref{eqn:hy_sys_def} evolve on hybrid time domains (HTDs) { \color{black}, denoted as \text{dom} ($\xi$)}, see \cite[Ch. 2]{goebel2012hybrid} or the Appendix in the Supplemental Material. 
We shall also utilize a class of hybrid inclusions with inputs of the form
\begin{subequations}\label{eqn:hy_sys_def_input}
\begin{align}
&(\xi,u)\in C\times\mathcal{U},~~~~~~\dot{\xi} \in F(\xi,u), \label{flowsHDS_input}\\
&(\xi,u)\in D\times\mathcal{U},~~~~ \xi^+ \in G(\xi,u),\label{jumpsHDS_input}
\end{align}
\end{subequations}
where $\mathcal{U}\subset\mathbb{R}^m$ is the set of admissible control inputs, and the remaining data of the system are similarly defined. For a complete mathematical description of HTDs and the notion of solution to hybrid inclusions and hybrid inclusions with inputs, we refer the reader to the Appendix {\color{black}(see also \cite{sanfelice2021hybrid})}. 

\vspace{0.1cm}
{\color{black} Lastly, we also define hybrid systems that allow for small perturbations or inflation in their data. Such systems can used to model small disturbances and study robustness properties.

\begin{defn}[Inflation of a hybrid system]
\label{defn:inflation} 
 The inflation of a hybrid system $\mathcal{H}$ is another hybrid system $\mathcal{H}_{\epsilon}$ given by
\begin{align}\label{eqn:inflation}
    x \in C_\epsilon,~~\Dot{x}\in F_\epsilon(x),~~~~~~
    x \in D_\epsilon,~~x^+\in G_\epsilon(x),  
\end{align}
where $\epsilon>0$, with data: 
\begin{align*}
    C_\epsilon &\coloneqq \{x \in \mathbb{R}^n \, : \, (x+\epsilon \mathbb{B}) \cap C \neq \emptyset\}, \\
    F_\epsilon(x) &\coloneqq \overline{\text{co}}~F_0((x+\epsilon \mathbb{B}) \cap C) + \epsilon \mathbb{B}, \\
    D_\epsilon &\coloneqq \{x \in \mathbb{R}^n \, : \, (x+\epsilon \mathbb{B}) \cap D \neq \emptyset\}, \\
    G_\epsilon(x) &\coloneqq G_0((x+ \epsilon \mathbb{B}) \cap D) + \epsilon \mathbb{B}.
\end{align*}   
\end{defn}
}

\vspace{0.1cm}
The following stability notions will be used in the paper:

\vspace{0.1cm}
\begin{defn}[Uniform Global Asymptotic Stability] 
    The set $\mathcal{A}$ is said to be uniformly globally asymptotically stable (UGAS) for system \eqref{eqn:hy_sys_def} if there exists $\beta \in \mathcal{KL}$ such that any solution $\xi$ satisfies the bound: 
    \begin{equation*}
        \lvert \xi(t,j) \rvert _{\mathcal{A}} \leq \beta(\lvert \xi(0,0) \rvert_\mathcal{A},t+j),  
    \end{equation*}
    for all $(t,j)\in\text{dom}(\xi)$. \QEDB 
\end{defn}

\vspace{1mm}
The following definition will be instrumental for the study of hybrid systems \eqref{eqn:hy_sys_def} that depend on a small parameter $\delta>0$.

\vspace{0.1cm}
{\color{black}
\begin{defn}[Semi-global Practical Asymptotic Stability]\label{defn:SGPAS} 
    The set $\mathcal{A}$ is said to be semi-globally practically asymptotically stable (SGPAS) for system \eqref{eqn:inflation} as $\epsilon \to 0$ if there exists $\beta \in \mathcal{KL}$ and for each $\mu,\nu >0$ there exists a $\epsilon^* >0$ such that, for each $\epsilon \in (0,\epsilon^*)$, any solution $\xi$ that satisfies $\lvert \xi(0,0) \rvert_{\mathcal{A}} \leq \mu$ also satisfies the bound
\begin{equation*}
    \lvert \xi(t,j) \rvert_{\mathcal{A}} \leq \beta(\lvert \xi(0,0) \rvert_{\mathcal{A}},t+j) + \nu,
\end{equation*}
for all $(t,j) \in \text{dom}(\xi)$. \QEDB 
\end{defn}
}

The following definition, which is fairly standard in the dynamical systems literature, will be instrumental in characterizing the stable set that arises in the switching systems considered in this paper.
\begin{defn}[Reachable set]\label{definitionReachable}
    Let $K \subset \mathbb{R}^n$ be a compact set. The (infinite horizon) reachable set from $K$ for the hybrid inclusion \eqref{eqn:hy_sys_def}, denoted as $R(K)$, is defined as
\begin{equation*}
    R(K) \coloneqq \{z \in \mathbb{R}^n  :  z=\xi(t,j), \xi \in \mathcal{S}(K), (t,j) \in \text{dom}(x)\}, 
\end{equation*}    
where $\mathcal{S}(K)$ denotes the set of solutions to system \eqref{eqn:hy_sys_def} with initial conditions in $K$. \QEDB 
\end{defn}

\vspace{0.1cm}
In this paper, some of the stability results will be defined with respect to \emph{$\Omega$-limit sets from sets}.  The following definition formalizes this concept.
\begin{defn}[$\Omega$-limit set]\label{defn:omega_set}
    Given a set $K$, the $\Omega$-limit set from $K$, denoted by $\Omega(K)$ is defined as 
\begin{align*}
    \Omega(K) \coloneqq \{&z \in \mathbb{R}^n \, :\, z = \lim_{i \to \infty} \xi_i(t_i,j_i),\\
    &\xi_i \in \mathcal{S}(K), (t_i,j_i) \in \text{dom}(\xi_i), \lim_{i \to \infty} t_i + j_i = \infty\}.
\end{align*}   
\end{defn}
\section{Problem Formulation}
\label{sec:prob_form}
In this section, we outline the model of the power transmission grid considered in this paper.

\subsection{Power Transmission System Model}
\label{subsecproblem}
We consider a power transmission grid with buses $\cN := \{1,...,N\}$ and lines $\cE \subset \cN \times \cN$. Let $\mathcal{D} \subset \cN$ be the set of buses where inverter-interfaced DEAs are connected, and let $\mathcal{G} \subset \cN$ be the set of buses where CGs are located.  We assume, without loss of generality, that $\cN = \cD\cup\cG$ and $\cD\cap\cG = \emptyset$. The set $\cI_{\ell} \subset \cE$ collects the lines connected to the bus $\ell$. We adopt a DC approximation of the power flows; transmission lines are lossless and  the reactance  of the line connecting bus $i$ and bus $j$ is denoted as $X_{ij}>0$. Next, we outline the models for both CGs and DEAs.

\subsubsection{CGs and frequency-responsive DEAs}
We start with the model for the CGs. Assuming that the exciter operates at a stable output such that the terminal voltage magnitude is constant, the following model for the CG $g \in \cG$ is widely adopted in the literature~\cite{kundur2007power}:
\begin{subequations} \label{eq:gen}
    \begin{align} 
    \dot\delta_g & = \omega_\mathrm{s}\Delta \omega_g,  \label{eq:gen_theta}\\
    M_{g} \Delta \dot \omega_g & =   P^{\mathrm{m}}_{g}  - D_{g} \Delta \omega_g - P_{\text{load},g} - \sum_{\ell \in \cI_g} P_{g \ell},\label{eq:gen_omega}\\  
    \tau_{g} \dot P^{\mathrm{m}}_{g} &= - P^{\mathrm{m}}_{g} + P_{g}^\mathrm{r} - K_{\mathrm{gov},g}  \Delta \omega_g, \label{eq:gen_p} 
    \end{align}
\end{subequations}
where $\delta_g$, $\Delta \omega_g$, $\omega_\mathrm s$, and $P^{\mathrm{m}}_{g}$ are the rotor electrical angle, the rotor speed deviation in per unit, the synchronous angular speed, and the turbine mechanical power, respectively.
Furthermore, $M_{g}$ is the constant of inertia, and $D_{g}$ models the equivalent load damping, which includes the damper windings. The dynamics of the turbine mechanical power are captured by a first-order turbine model \cite{dorfler2023control}, where $K_{\mathrm{gov},g}$ is the governor gain, modeling the inverse of the speed-droop regulation constant, $\tau_g$ is the turbine time constant, and $P^{\mathrm{r}}_{g}$ denotes the reference-power setting computed from a higher layer control (e.g., a secondary controller). Finally, $P_{\text{load},g}$ is the real load at bus $g$, and $P_{g \ell} = \frac{\delta_g - \delta_\ell}{X_{g \ell}}$ is the power flow from bus $g$ to $\ell$.


On the other hand, for each frequency-responsive DEA $d \in \mathcal D$, we consider the following dynamics~\cite{guggilam2018optimizing,wang2024distributed,jung2024inferences}:
\begin{subequations} \label{eq:DER}
    \begin{align} 
    \dot \delta_d &= \omega_\mathrm s \Delta \omega_d \label{eq:der_theta}, \\
    M_{d} \Delta \dot \omega_d &=   - D_{d} \Delta \omega_d + P_{\text{DEA},d}^\text{r} - P_{\text{load},d} - \sum_{\ell \in \cI_d}P_{d \ell}, \label{eq:der_omega}  
    \end{align}
\end{subequations}
where $D_{d}$ models the frequency response of the DEA and
$M_{d}$ determines the (virtual) inertial response, and $P_{\text{DEA},d}^\text{r}$ and $P_{\text{load},d}$ denote the reference active power of the DEA and the total load at node $d$, respectively. We note that $D_d$ and $M_d$ do not represent mechanical parameters as in~\eqref{eq:gen}. Instead, for DEAs, these are \emph{digital} parameters that may be tuned to obtain a desired response \cite{mavalizadeh2023improving, cheng2020smart}.

Before proceeding, we note that uncontrollable frequency-sensitive loads are not modeled
to simplify notation, although they can be straightforwardly incorporated into~\eqref{eq:DER}. 
The dynamics \eqref{eq:gen} and \eqref{eq:DER} lead to a linear time-invariant (LTI) system of the form 
\begin{equation}
    \label{eqn:ss_full}
    \Dot{x}_\text{f} = \bar{A}x_\text{f} + \bar{B}u_\text{f}, 
\end{equation}
with $x_\text{f} \coloneqq ( \delta_1 , \ldots, \delta_N, \Delta \omega_1, \ldots, \Delta \omega_N, P_1^\text{m}, \ldots, P_{|\mathcal{G}|}^\text{m}) \in \mathbb{R}^{3\lvert \mathcal{G} \rvert + 2 \lvert \mathcal{D} \rvert}$,  and where the vector $u_\text{f} \in \mathbb{Z}_{\geq 0}$ collects the loads $P_{\text{load},g}$ for $g \in \mathcal{G}$ and the net loads $ P_{\text{load},d} - P_{\text{DEA},d}^\text{r}$ at the nodes $d \in \mathcal{D}$. The   matrices $\bar{A}$ and $\bar{B}$ can be readily constructed from \eqref{eq:gen} and \eqref{eq:DER} (their explicit expression is omitted due to space limitations; see e.g.,~\cite{jung2024inferences,colombino2019online}).

Model reduction and aggregation have been widely adopted in the power systems literature to simplify analysis and simulations {\color{black}(see \cite{ajala2022model})}. In this paper, to reduce the number of state variables and to facilitate the analysis, we adopt the model of~\cite{guggilam2018optimizing,jung2024inferences}. In particular,  we assume that $\Delta \omega$ is (approximately) the same for all nodes; this assumption is reasonable for networks where electrical distances are negligible and all the buses have the same frequency even during transients; see, e.g., \cite{anderson1990low}. 

{\color{black} Moreover, we deliberately neglect the fast dynamics of the inverter inner control loops—specifically, the voltage and current regulators—as our focus is on the electromechanical behavior of the aggregated system under switching net loads (see, for example, the considerations make by~\cite{dorfler2023control}). The inner loops typically operate at bandwidths in the kilohertz range, which are several orders of magnitude faster than the outer power and frequency control loops (see, for example, \cite{Yuantime} and \cite{he2020transient}). As a result, their dynamics can be regarded as quasi-static on the electromechanical timescale. This separation of time scales, well supported by the theory of singular perturbations (\cite{kokotovic1999singular,chow1982time}), is a widely accepted assumption in power-system studies concerned with frequency stability and primary control.}

{\color{black} For example, in \cite[Section 2.2]{samanta2025stability}, the authors consider small-signal models of voltage-source converters with output LC filters and specifically ignore the fast inner loop dynamics stating that they are significantly faster than primary control. A similar argument is stated in \cite{bakhshizadeh2023improved} wherein reduced-order models of wind turbines are developed to study instabilities induced by phase locked loops. For a more technical treatment, we refer the reader to \cite{schiffer2016survey}, see in particular Assumption 5.1.}

Then, for a lossless system, we have that $\sum_{g \in \cG}\sum_{\ell \in \cI_g}P_{g \ell} + \sum_{d \in \cD}\sum_{\ell \in \cI_d}P_{d \ell}=0$. Using~\eqref{eq:gen} and~\eqref{eq:DER} we obtain
\begin{equation}
M_{\mathrm{eff}} \Delta \dot \omega = \sum_{g \in \cG}  P^{\mathrm{m}}_{g}  -D_\mathrm{net} \Delta \omega - P_\mathrm{load},\label{eq:commonfreq}
\end{equation}
where  $
    P_{\mathrm{load}} \coloneqq -\sum_{g \in \cG} P_{\mathrm{load},g} - \sum_{d \in \cD} (P_{\mathrm{load},d} - P_{\mathrm{DEA},d})$ is the aggregate electrical load across the transmission system (demand minus the total  generation from DEAs), and where the \emph{effective inertia constant} $M_\mathrm{eff}$ and the \emph{net damping constant} $D_\mathrm{net}$ are defined, respectively, as:
    \begin{align}
    \label{eq:EffParams}
        \hspace{-.25cm} M_{\mathrm{eff}} \coloneqq \sum_{g \in \cG} M_{g} + \sum_{d\in \cD} M_{d}, ~ 
        D_{\mathrm{net}} \coloneqq \sum_{g \in \cG} D_{g} + \sum_{d\in \cD} D_{d}. 
    \end{align}    
Furthermore, from \eqref{eq:gen_p} we have:
\begin{equation}
\diag(\tau) \dot P_\cG^{\mathrm{m}} = - P_\cG^{\mathrm{m}} + P_\cG^\mathrm{r} - K_{\mathrm{gov},\cG} \Delta \omega, \label{eq:commonfreqP}
\end{equation}
where $\tau$ is a vector collecting $\{\tau_i\}_{i \in \cG}$ and $K_{\mathrm{gov},\cG}$ is a vector collecting $\{K_{\mathrm{gov},i}\}_{i \in \cG}$.

We will utilize this aggregated model in the subsequent analysis. However, in the numerical results that will be presented in Section~\ref{sec:sims}, we will test both the aggregate model and the model~\eqref{eqn:ss_full} in order to validate our analysis. 
\begin{figure}
    \centering
    \includegraphics[width=\linewidth]{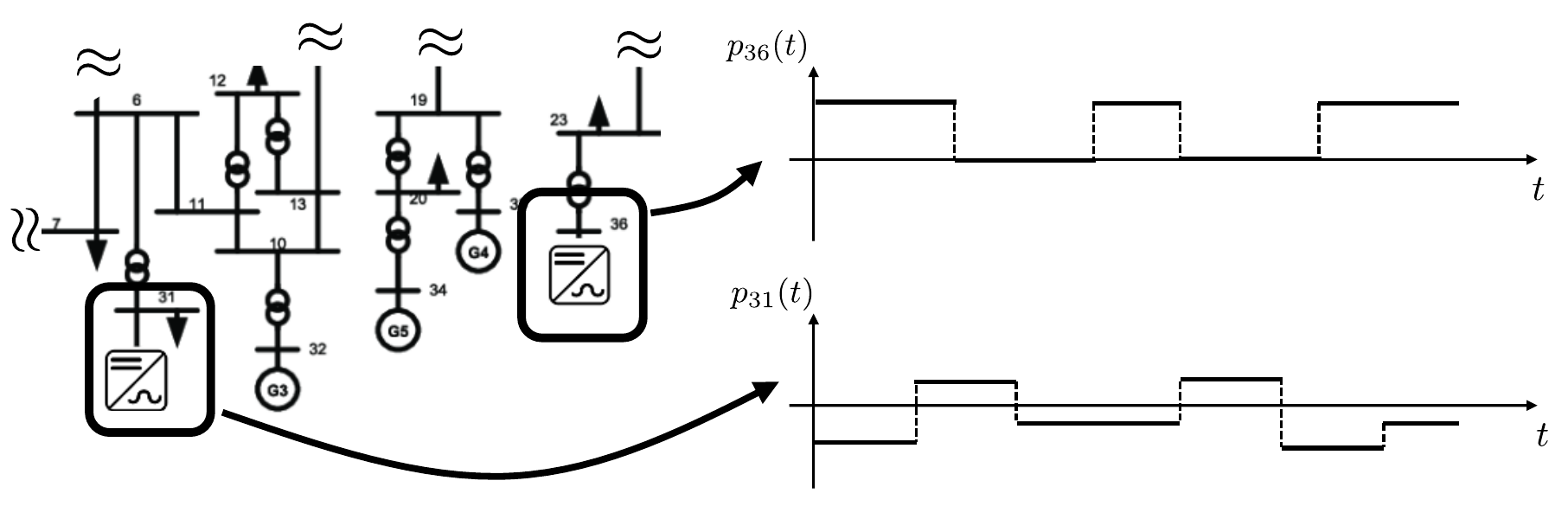}
    \caption{Representation of transmission power system under switching net loads.}
    \label{fig:bd}
    \vspace{-.3cm}
\end{figure}
%

\subsubsection{Secondary controller}
The reference powers of the generators $P_{g}^\mathrm{r}$, $g \in \mathcal{G}$ are computed by a secondary frequency controller, whose objective is to steer the frequency deviation to zero~\cite{kundur2007power,wood2013power}. We consider the model of~\cite[Ch.~9]{wood2013power}, and let $z \in \mathbb{R}$ be the   state of the secondary controller. The dynamics of this controller are  given by:
\begin{subequations} \label{eqn: secondary_ctrl}
\begin{align}
\tau_z \dot{z} &= -z + \beta \Delta \omega + \one^\top P_\cG^\mathrm m,\\
P_\cG^\mathrm r &= P_\cG^* + \zeta (z-\one^\top P_\cG^*),
\end{align}
\end{subequations}
where $\zeta \in \mathbb{R}^{|\cG|}_{\geq 0}$ are the participation factors (i.e., $\zeta_i \in [0,1]$ and $\one^\top \zeta = 1$), $\beta \in \mathbb{R}_{< 0}$ is a tunable gain,  $P^{\mathrm{m}}_{\cG} := [\{P^{\mathrm{m}}_{i}\}_{i \in \cG}]^\top$ collects the mechanical powers of every generator $g\in \cG$, and where  $P_\cG^* \in \mathbb{R}^{|\cG|}$ are operating points for the CGs that are computed via tertiary control (e.g., economic dispatch or DC optimal power flow).  
To ensure a sufficient time-scale separation between the primary and secondary frequency control, we have that  $\tau_z > \tau_g, \,  \forall ~g \in \cG$. 

\vspace{0.1cm}
\subsubsection{State-space Model}

Combining~\eqref{eq:commonfreq},~\eqref{eq:commonfreqP} and~\eqref{eqn: secondary_ctrl} yields the state-space model of the power transmission system:
\begin{equation}\label{eqn:switched_loads}
        \dot{x} =A x + B u,
\end{equation}    
where the state vector is $x = (\Delta \omega, P_\cG^\mathrm m ,z)\in \mathbb{R}^{(2+|\cG|)}$, the input is defined as  $u = (P_\mathrm{load}, P_\cG^*)\in \mathbb{R}^{(1+|\cG|)}$, and the matrices $A \in \mathbb{R}^{(2+|\cG|) \times (2+|\cG|)}$ and $B \in \mathbb{R}^{(2+|\cG|) \times (1 +|\cG|)}$ are given by: 
\begin{align}
\label{eqn: ss_new}
A &= \left[
 \begin{array}{c c c}
-D_{\mathrm{net}}M_{\mathrm{eff}}^{-1}  & M_{\mathrm{eff}}^{-1} \one^\top & 0\\
A_\tau K_{\mathrm{gov},\cG} & A_\tau & -A_\tau \zeta\\
\tau_z^{-1} \psi & \tau_z^{-1} \one^\top & -\tau_z^{-1} \\
 \end{array} 
 \right],\\
B &=   \left[
\begin{array}{c c}
-M_{\mathrm{eff}}^{-1}  &  \zero^\top\\
0 & -A_\tau \left(\mathbb{I}-\zeta\one^\top\right) \\
0 & \zero^\top \\
 \end{array} 
 \right], \nonumber
\end{align}
where $A_\tau := -\diag(\tau_1, \ldots, \tau_{|\mathcal{G}|})^{-1}$. One can verify that, for typical power systems setups, the matrix $A$ is Hurwitz stable. {\color{black} Additionally, even if some of the (virtual) inertia and damping parameters are set to zero, the matrix $A$ remains Hurwitz as long as damping and inertia are provided in the aggregate. That is, $D_{\text{net}}, M_{\text{eff}} > 0$ even though $D_d = M_d = 0$ for some $d \in \mathcal{D}$. The case of $D_{\text{net}} = 0$ leads to the matrix $A$ being only marginally stable. This could result in sustained oscillations in certain system states. As such, we consider \emph{strictly} positive virtual damping and inertia to be consistent with physical reality (existing power systems always have non-zero mechanical and electrical losses) and to be compliant with interconnection standards of DEAs to power system interfaces like \cite{IEEE-Std-1547-2018}. Moreover, a non-zero damping term proportional to the system frequency is always included in standard references on power systems modeling (see \cite{kundur2007power, dorfler2023control}).}

\subsubsection{Switching net loads}
\label{subsec:switched}
Recall that $P_{\mathrm{load}} = -\sum_{g \in \cG} P_{\mathrm{load},g} - \sum_{d \in \cD} (P_{\mathrm{load},d} - P_{\mathrm{DEA},d})$ is the aggregate electrical load across the transmission system (demand minus the total  generation from DEAs). We decouple  $P_{\text{load}}$ as $P_{\text{load}} = P_{\text{switching}} + P'_{\text{load}}$, where $P'_{\text{load}}$ are aggregate loads that vary continuously over time, while $P_{\text{switching}}$ models aggregate net loads that may exhibit a switching behavior.  Consequently, in \eqref{eqn:switched_loads}, the input is given by $u = u _{\text{switching}} + \tilde{u}$. Hereafter, we assume that $u_{\text{switching}}$ takes values from a finite set $\{u_q\}_{q\in\mathcal{Q}}$, indexed by ${ \color{blue} \mathcal{Q}}$. 

The presence of the component $P_{\text{switching}}$ is motivated by dynamic operating conditions
that may result from protection scheme activations, external cyber-attacks, or due to the integration of dynamic loads,
such as data centers: 

\noindent $\bullet$ Cyber-attacks may change the power setpoints $P_{\mathrm{DEA},d}^\text{r}$ of inverter-interfaced DEAs intermittently. Denial-of-service attacks can also periodically shut down the inverters (see, e.g., \cite{kang2015investigating}, \cite{zografopoulos2021detection}, and \cite{barua2020hall}). 

\noindent $\bullet$ When overvoltages occur in feeders, protection schemes may deactivate inverters, leading to changes in the net loading at the substation. For instance, the standard CENELEC EN50549-2 \cite{standard} defines three statuses for the  DEAs: running, idling, and disconnected. In the idling or disconnected states, no active power is supplied by the DEA. After a specified duration, the inverter reconnects, which can result in periodic fluctuations in the power levels at the substation.

\noindent $\bullet$ Dynamic loads, particularly data centers, can display switching behavior. This is especially true for loads like cryptocurrency mining, where power consumption can have a given number of power ``levels'' \cite{wheeler2018power}. 

{\color{black}
\noindent $\bullet$ Changes in the optimal power flow (OPF) schedule translates to updating the generator setpoints $P^*_{\mathcal{G}}$. Typically, these variations occur at the tertiary time-scale and can be treated as slowly-varying changes to the input term $\Tilde{u}$ in Equation \eqref{eqn:switched_loads}. In contrast, if such changes are abrupt (cyber-attack contingencies or protection-triggered re-dispatch), the updated generator setpoints can be absorbed into the term $u_{\text{switching}}$ leading to additional modes in the set $\mathcal{Q}$. 
}

\vspace{0.2 cm}
We note that, in the first two cases, the values of $D_{\mathrm{net}}$ and $M_{\mathrm{eff}}$ change and, accordingly, the matrices $A$ and $B$ also depend on $q$. Therefore, hereafter we will use the notation $A_q$ and $B_q$, $q \in \mathcal{Q}$ to denote the specific values of the matrices in the mode $q$. The dynamics of the resulting system are then captured by a system that switches between a finite number of vector fields, or modes. For a given $q \in \mathcal{Q}$,~\eqref{eqn:switched_loads} then becomes:  
\begin{align}
    \label{eqn: ss_switched}
    \Dot{x} &= A_q x + \underbrace{B_q u_q}_{:=b_q} + B_q \tilde{u}  \, .        
\end{align}    

\noindent 

To further clarify our setup, take as an example a system where DEAs are located at buses 31 and 36 as in see Figure \ref{fig:bd}; consider the case where DEAs undergo an intermittent denial-of-service attack. 
Then, in this example, $u_q = ( P_{\mathrm{DEA},31} + P_{\mathrm{DEA},36}, 0)$ and  $\tilde{u} = (-\sum_{n \in \cN} P_{\mathrm{load},g} + \sum_{d \in \cD \backslash \{31, 36\}} P_{\mathrm{DEA},d}, P_\cG^*)$. 
%
{
\color{black}
\begin{com}
    Nodes with only conventional loads are governed by the same model as in~\eqref{eq:DER} with $D_{d} = 0$ and $M_d = 0$. However, setting the damping and inertia parameters to zero, reduces \eqref{eq:der_omega} to an algebraic constraint that models power balance and would lead to a system of differential \emph{algebraic} equations (DAEs) as opposed to the system of first order linear differential equations (with switching vector fields and multiple equilibria) that we consider in Equation \eqref{eqn: ss_switched}. While several works address the stability analysis of switched DAEs (\cite{zhou2013stability, liberzon2012switched, debeljkovic2011stability}), they do so only with respect to a \emph{single} equilibrium point. The case of multiple equilibria is still largely unexplored. Therefore, extending our current theoretical treatment to buses with only conventional loads is a subject for future work.
\end{com}}
\vspace{-0.3cm}
{\color{black}
\begin{com}
    Topology changes of the transmission network - line outages and breaker switching, can also be captured by the model in \eqref{eqn: ss_switched} provided the electrical distances among buses, post-change, still remain negligible (see \cite{cotilla2013multi} and \cite[Chapter 3]{chow2013power}). Such reconfigurations typically result in parametric changes in damping, inertia, and the participation factors resulting in an additional system mode via the matrices $A_q$ and $B_q$. Variations that result in weakly coupled areas necessitate a multi-area extension of our current model via the inclusion of a susceptance matrix and is a subject of future work.
\end{com}
}
\vspace{-0.2cm}
\section{Main Results}
\label{sec:main}
In this section, we formulate and prove our main results concerning the stability properties of the power transmission system \eqref{eqn: ss_switched} under intermittent loads. 
Throughout the current section, we impose the following assumption.
\begin{assump}\label{asmp:hurwitz_Aq}
    For all $q\in\mathcal{Q}$, the matrix $A_q$ in \eqref{eqn: ss_switched} is Hurwitz. \QEDB
\end{assump}
We follow two distinct approaches in assessing the stability properties of system \eqref{eqn: ss_switched}. In the first approach, which appears in Subsection \ref{subsec:approach_1} below, we consider the case where the intermittent load is \textit{slowly varying}, and we establish a \emph{practical stability} result for the switched system with respect to a suitable compact set. In the second approach, which appears in Subsection \ref{subsec:approach_2}, we allow for arbitrarily fast changing loads and obtain a conservative input-to-state stability result with respect to a larger set.

\subsection{Slowly Varying Loads}\label{subsec:approach_1}
To study the stability properties of system \eqref{eqn: ss_switched} under slowly varying load and under sufficiently slow switching, we leverage the framework of Hybrid Dynamical Systems (HDS) \cite{goebel2012hybrid}.
Specifically, we introduce the variable $\tau \in \mathbb{R}_{\geq 0}$ to serve as the state of a timer responsible for triggering mode transitions. Then, we model the switching signal in \eqref{eqn: ss_switched} as a solution $q$ of the following hybrid automaton with set-valued dynamics:
\begin{subequations}\label{eq:swtching_automata}
    \begin{align}
       &(q,\tau)\in\mathcal{Q}\times[0,1],~~~~~~(\dot{q},\dot{\tau})\in\{0\}\times[0,\delta_1] \\
        &(q,\tau)\in \mathcal{Q}\times\{1\},~~~~(q^+,\tau^+)\in (\mathcal{Q}\backslash\{q\})\times \{\tau-1\}.
    \end{align}
\end{subequations}
As shown in \cite[Sec. 2.4]{goebel2012hybrid}, every solution $q$ of \eqref{eq:swtching_automata} satisfies a dwell-time constraint with dwell time given by $\delta_1^{-1}$, i.e., any two switches of $q$ are separated by at least $\delta_1^{-1}$ Moreover, every switching signal with such a dwell time can be generated via \eqref{eq:swtching_automata}.

Similarly, we model the time-varying load $\tilde{u}$ as a signal generated by an \textit{exogenous} dynamical system of the form
\begin{align}\label{eq:continuous_load}
   \tilde{u}\in\mathcal{U},~~~~  \dot{\tilde{u}} \in \delta_2 \Pi(\tilde{u}),
\end{align}
where $\mathcal{U}$ is a compact set, $\delta_2\geq 0$ is a parameter that characterizes the rate of change of $\tilde{u}$, and $\Pi$ is any {\color{black} multi-valued } function that satisfies the following assumption:

{\color{blue}\vspace{0.1cm}
\begin{assump}
    \label{assump:exo}
    The set $\mathcal U$ is compact, and the set-valued map $\Pi:\mathcal U\rightrightarrows\mathbb R^{m}$ is OSC, LB, and
convex valued. Moreover, $\mathcal U$ is strongly forward invariant
for the dynamics $\dot{\tilde u}\in\delta_2\Pi(\tilde u)$. \QEDB 
\end{assump}}
%

\begin{rem}
    Assumption \ref{assump:exo} guarantees that the input $\tilde{u}$ is uniformly bounded, although the uniform upper bound can be arbitrary since it depends on the size of the compact set $\mathcal{U}$. On the other hand, the continuous-time evolution defined by the differential inclusion \eqref{eq:continuous_load} is fairly general and admits a rich class of signals $\tilde{u}$ of interest in power systems. \QEDB 
\end{rem}

\vspace{0.1cm}
Next, we introduce the shifted coordinates
\begin{align}\label{eq:shfted_coords}
    y = x + A_q^{-1}B_q\tilde{u},
\end{align}
which represent the deviation of the state $x$ from the point $-A_q^{-1}B_q\tilde{u}$, which is the equilibrium point of \eqref{eqn: ss_switched} when $b_q=0$ {\color{black} and $\Tilde{u}$ is constant}. A direct computation shows that the continuous-time evolution of the state $y$ is governed by the differential inclusion
\begin{align}\label{eq:flow_shifted_coords}
    \dot{y}\in \{A_q y + b_q\} + \delta_2 A_q^{-1}B_q \Pi(\tilde{u}).
\end{align}
To compute the discrete-time dynamics of the state $y$, we note that whenever the mode $q$ switches to $q^+\in\mathcal{Q}\backslash\{q\}$, the jumps of $y$ satisfy
\begin{align*}
    y^+= y + \left(A_{q^+}^{-1}B_{q^+}-A_{q}^{-1}B_{q}\right)\tilde{u}.
\end{align*}
By defining the parameter $\delta_3\geq 0$ by
\begin{align*}
    \delta_3:=\max_{q\in\mathcal{Q}} \max_{q^+\in\mathcal{Q}\backslash\{q\}} |A_{q^+}^{-1}B_{q^+}-A_{q}^{-1}B_{q}|~\cdot ~\max_{\tilde{u}\in\mathcal{U}}|\tilde{u}|,
\end{align*}
we observe that $y^+$ satisfies the following difference inclusion: 
\begin{align}\label{eq:jump_shifted_coords}
    y^+ \in \{y\} + \delta_3\mathbb{B}.
\end{align}
Therefore, by introducing the parameter $\delta:=\max\{\delta_1,\delta_2,\delta_3\}$ and combining \eqref{eq:swtching_automata}-\eqref{eq:jump_shifted_coords}, we can obtain the overall hybrid dynamics with state $\xi \coloneqq (y, q, \tilde{u}, \tau)$ and data $\mathcal{H}_{\delta}=(C,F_\delta,D,G_\delta)$, where the flow set $C$ and the jump set $D$ are given by
\begin{subequations}
\label{eqn:hybrid_system}
\begin{align}
    \label{eqn:flow_set}
    &C \coloneqq \mathbb{R}^{\lvert \mathcal{G} \rvert + 2} \times \mathcal{Q} \times \mathcal{U} \times [0,1],\\
    \label{eqn:jump_set}
    &D \coloneqq \mathbb{R}^{\lvert \mathcal{G} \rvert + 2} \times \mathcal{Q} \times \mathcal{U}  \times \{1\},
\end{align}
and the flow map $F_\delta$ and the jump map $G_\delta$ are
\begin{align}
    \hspace{-5cm}F_\delta(\xi) &\coloneqq \Big(\{(A_qy + b_q)\}+\delta_2 A_q^{-1} B_q\Pi(\tilde{u})\Big) \notag\\
    & ~~~~~~~~~~~~~~~~~~~~~~~~\times \{0\}\times \delta_2\Pi(\tilde{u}) \times
        [0, \delta_1],\label{eqn:flow_map}\\
    G_\delta(\xi) &\coloneqq \Big(\{y\} + \delta_3 \mathbb{B}\Big)\times
        (\mathcal Q\backslash\{q\})\times \{\tilde{u}\}\times\{\tau -1\}.\label{eqn:jump_map}
\end{align}
\end{subequations}
With the model \eqref{eqn:hybrid_system} at hand, we can now present the first main result of the paper.

\vspace{0.1cm}
\begin{thm}
    \label{thm:stability_result}
    Let Assumptions \ref{asmp:hurwitz_Aq} and \ref{assump:exo} hold. Define the set 
{\color{blue}\begin{align}\label{eqn:omega_limit_set_1}
\mathcal{A}&:=\bigcup_{q\in\mathcal{Q}}\Big[
\{y_q^*\}\times\{q\}\times\mathcal{U}\times[0,1]
\notag\\
&\quad\cup
\bigcup_{r\in\mathcal{Q}\setminus\{q\}}
\overline{\{\Theta_t^r(y_q^*):t\geq0\}}
\times\{r\}\times\mathcal{U}\times\{0\}
\Big].
\end{align}}
where the map ${\color{black}\Theta^q_t}:\mathbb{R}^{|\mathcal{G}|+2}\rightarrow \mathbb{R}^{|\mathcal{G}|+2}$ is given by
\begin{align}\label{eq:flow_individual_linear_mode}
    \Theta_t^{q}(y):= y_q^* + \exp(A_q t) (y - y_q^*),~~~~~y_q^*:=- A_q^{-1} b_q.
\end{align}
Then, the set $\mathcal{A}$ is compact and semi-globally practically asymptotically stable as $\delta\to0^+$ for $\mathcal{H}_{\delta}$. \QEDB 
\end{thm}

\vspace{0.1cm}
The result of Theorem \ref{thm:stability_result}, {\color{black}which relies on stability tools for switching systems with distinct equilibrium \cite{baradaran2020omega}, as well as a sharp characterization of the so-called $\Omega$-limit set of the resulting system (presented below in Proposition 1),} states the stability properties of the hybrid system \eqref{eqn:hybrid_system} with respect to the non-trivial compact set $\mathcal{A}$, and under sufficiently slow switching and load variation, and sufficiently small mode-dependent mismatch induced by the coordinate transformation, i.e., small values of $\delta$. {\color{blue} Define the set $\mathcal{Y} \coloneqq \text{proj}_y(\mathcal{A})$, i.e., the projection of the set $\mathcal{A}$ onto the $y$-component of the system $\mathcal{H}_\delta$, and note that
\begin{equation}
    \mathcal{Y}
    := \bigcup_{q\in\mathcal{Q}} \overline{\mathcal{Y}_q},
    \qquad
    \mathcal{Y}_q
:=
\bigcup_{r\in\mathcal{Q}\setminus\{q\}}
\left\{
    \Theta_t^r(y_q^\ast):t\geq0
\right\}.
    \label{eq:Y}
\end{equation}
This set describes} the limiting behavior of the state variables in \eqref{eqn:switched_loads} in terms of the equilibrium points of the individual subsystems of \eqref{eqn: ss_switched} for sufficiently slowly occurring DEA-related events. We remark that the structure of the set $\mathcal{A}$, completely decouples the state $y$ from the exogenous input $\tilde{u}$. {\color{blue} It can be inferred from Equation \eqref{eq:Y}} that the shape of the set $\mathcal{Y}$ is independent of the exogenous input $\tilde{u}$. Such a decoupling is entirely due to the change of coordinate \eqref{eq:shfted_coords} and is not possible in the original coordinates. {\color{blue}Moreover, SGPAS of \(\mathcal A\) implies practical convergence of \(y\) to \(\mathcal Y\).}

\vspace{0.1cm}

{\color{black} Additionally, the existence of the set $\mathcal{A}$ as a compact, practically-stable, attractor does \emph{not} depend on the knowledge of system parameters, namely the matrices $A_q, B_q, \, \forall \in \mathcal{Q}$. Theorem \ref{thm:stability_result}, primarily, establishes a qualitative result regarding the long-term behavior of the state trajectories of system $\mathcal{H}_\delta$ for sufficiently small values of $\delta$ (see also Remark \ref{rem:stab_res_rem}) and arbitrary values of parameters like the effective damping and inertia. Of course, for a particular instance of $\mathcal{H}_\delta$, with known parameters, Theorem \ref{thm:stability_result} also provides a computationally-tractable method for a quantitative estimate of $\mathcal{A}$ via Equations \eqref{eqn:omega_limit_set_1}{\color{blue}, or \eqref{eq:Y} (for the physical state).}}  

\vspace{0.1cm}
Before discussing in more detail the structure of the set $\mathcal{A}$ and the proof of Theorem \ref{thm:stability_result}, we present a stylized numerical example that provides a graphical visualization of the set $\mathcal{Y}$ for a system of the form \eqref{eqn: ss_switched} with two modes, each mode generating a stable equilibrium point.

\begin{figure*}[t!]
    \centering
    \includegraphics[width=0.45\linewidth]{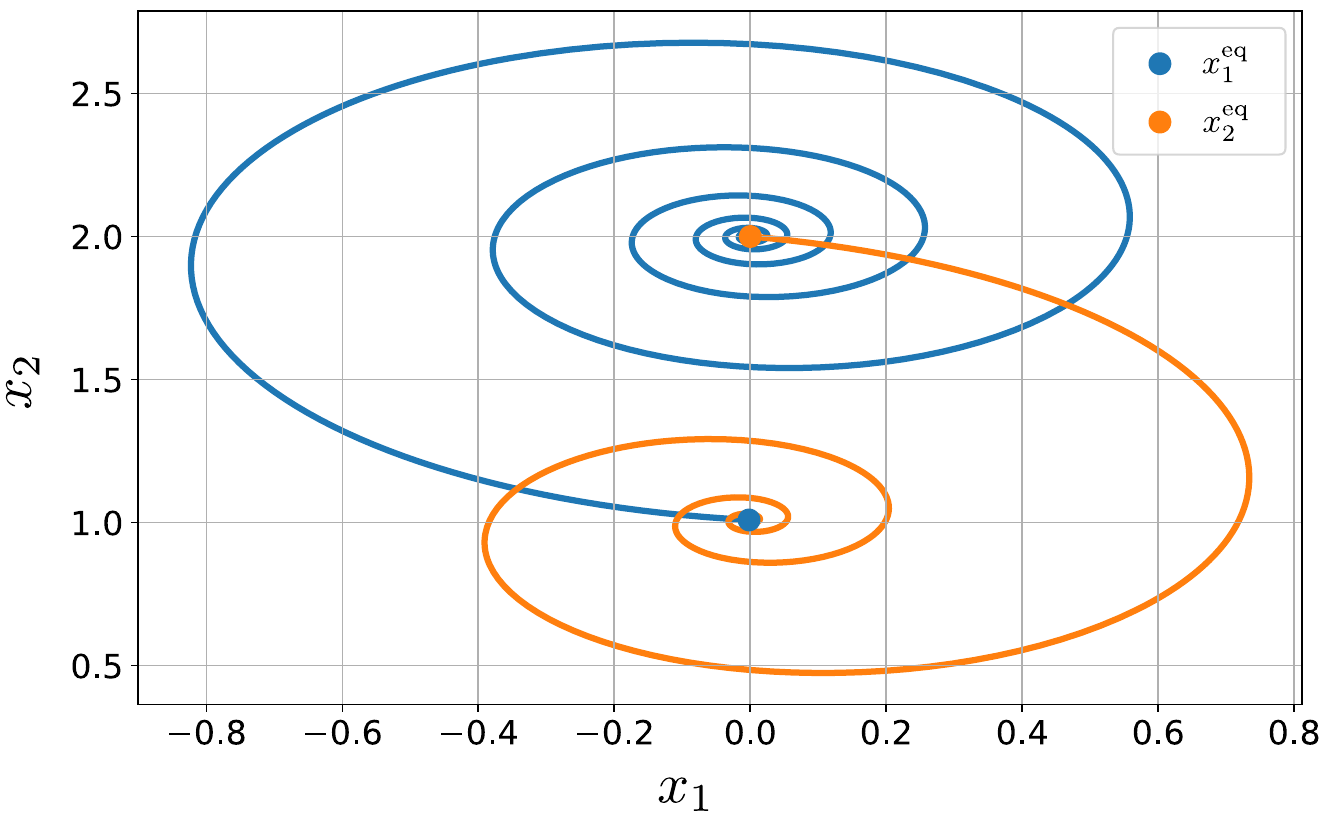}\includegraphics[width=0.45\linewidth]{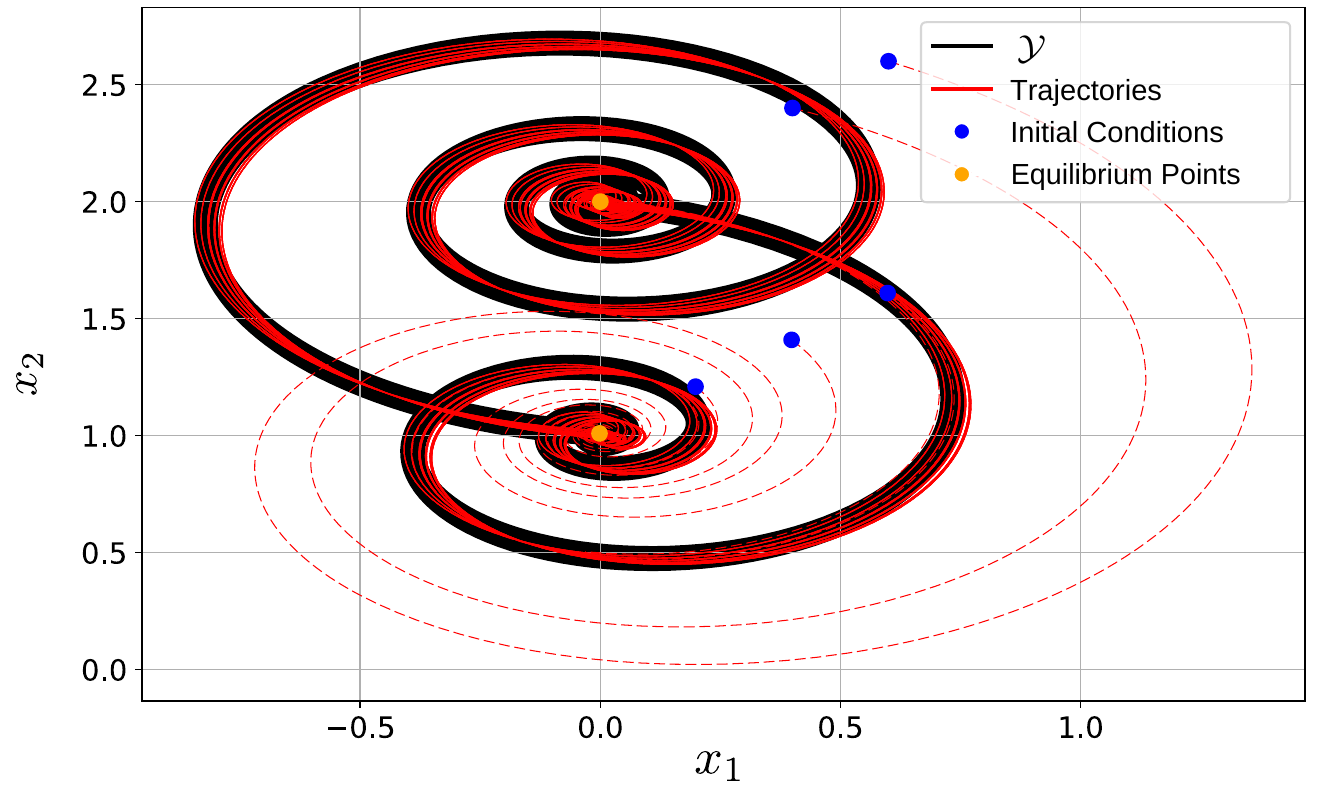}
    \caption{(Left) The set $\mathcal{Y}$ corresponding to the data in Example \ref{ex:nominal_toy}. (Right) Trajectories of the state $(\Delta \omega, P_m)$ with the data defined in Example \ref{ex:nominal_toy} under sufficiently small switching (red), overlaid on the set $\mathcal{Y}$ (thick black). As predicted by Theorem \ref{thm:stability_result}, the trajectories converge to an arbitrarily small neighborhood around $\mathcal{Y}$. }
    \label{fig:nominal_toy}
\end{figure*}

\vspace{0.1cm}
\begin{ex}
\label{ex:nominal_toy}
Consider a simplified power transmission system with one synchronous generator and a collection of DEAs, with the state vector $(x_1, ~ x_2) = (\Delta \omega , ~ P_m)$, modeled as $\Dot{x} = A_qx + b_q$, $q \in \{1,2\}$, where
\begin{align*}
    &A_1 = \begin{bmatrix}
    -0.6 & 2.98 \\
     -2.98 &  -0.6
\end{bmatrix}, ~~b_1 = \begin{bmatrix}
    -2.98 \\ 0.6
\end{bmatrix},\\
&A_2 = \begin{bmatrix}
    -0.4 & 3.24 \\
    -3.24 & -0.4 
\end{bmatrix}, ~~b_2 = \begin{bmatrix}
    -6.48 \\ 0.8
\end{bmatrix}. 
\end{align*}

In this example, we resort to the droop controller for frequency regulation, and for simplicity, we take $\mathcal{U}=\{0\}$ and $\Pi(u)=\{0\}$. As can be verified, both $A_1$ and $A_2$ are Hurwitz and, therefore, Assumption \ref{asmp:hurwitz_Aq} is satisfied. In particular, the two equilibrium points $x^1_{eq} = (0, 1)$ and  $x^2_{eq} = (0,2)$ are stable.
 Figure \ref{fig:nominal_toy} illustrates the shape of the set $\mathcal{Y}$. We simulate the hybrid system $\mathcal{H}_\delta$ for a sufficiently small $\delta$. Figure \ref{fig:nominal_toy} shows the asymptotic behavior of this system from five different initial conditions. It can be seen that all five trajectories closely approximate the set $\mathcal{Y}$. The frequency deviation remains close to its equilibrium point of zero in any of the two modes and is perturbed from zero precisely when the switching occurs. After the switching, the deviation is driven to zero. Similarly, the mechanical power oscillates between the two values of $1$ and $2$ (p.u.). Both these patterns of oscillations can be seen in Figure \ref{fig:nominal_toy} showcasing that in the presence of multiple equilibria and sufficiently slow switching, the {\color{blue} trajectories of the state $(\Delta \omega, P_m)$ remain arbitrarily close to the set $\mathcal{Y}$.}  \hfill $\Box$
\end{ex}

\vspace{0.1cm}
\begin{rem}
    \label{rem:stab_res_rem}
    {\color{black}
    According to Theorem \ref{thm:stability_result}, the trajectories of the power transmission system under slowly time-varying loads in addition to switching loads will, in the limit, closely approximate the set $\mathcal{A}$ from any initial configuration. The set $\mathcal{A}$ is not only (practically) attractive, but also stable. Namely, if the trajectories are initialized close to $\mathcal{A}$, they will remain close to $\mathcal{A}$ under slow switching of loads.}  \hfill \QEDB
\end{rem}

\vspace{0.1cm}
\begin{rem}
\label{rem:prop_1}
In the proof of Theorem \ref{thm:stability_result}, which follows shortly, we show that the set $\mathcal{A}$ is constructed as the $\Omega$-limit set of a certain ``nominal" hybrid system obtained in the limit $\delta\searrow 0$. 
By relying on stability results for switching and hybrid systems \cite{baradaran2020omega}, we show that the set $\mathcal{A}$ is uniformly globally asymptotically stable for $\mathcal{H}_0$. However, the structure of the set $\mathcal{A}$ in \eqref{eqn:omega_limit_set_1} does not immediately follow from the definition of the $\Omega$-limit set or the results in \cite{baradaran2020omega}, which are derived for a more general class of switching systems satisfying average dwell-time constraints. Instead, Theorem \ref{thm:stability_result} relies on Proposition \ref{prop:desc} (presented below), which shows that under the more typical dwell-time constraints imposed on switching loads in power transmission systems, the $\Omega$-limit set admits a sharper representation given by \eqref{eqn:omega_limit_set_1}. This representation relies on the mapping $\Theta_t^{q}(\cdot)$, which enables a simpler computational framework for the study of the asymptotic behavior of the system. More precisely, to compute $\mathcal{A}$, we can select an equilibrium point and compute the solution maps under the remaining subsystems starting from said point. Such a computation can be readily done using the variation of constants formula for linear systems. {\color{blue}Repeating this procedure for all the equilibrium points and taking their union with the mode and timer coordinates
specified in \eqref{eqn:omega_limit_set_1} will result in an explicit representation of $\mathcal{A}$.} Thus, Theorem 1 does away with having to perform simulations on shrinking sets of initial conditions and computing common points.
\hfill $\Box$
\end{rem}

\vspace{-0.2cm}

We are now ready to provide the proof of Theorem \ref{thm:stability_result} and begin by introducing a ``nominal'' version of the hybrid system $\mathcal{H}_{\delta}$ defined in \eqref{eqn:hybrid_system}, denoted by $\mathcal{H}_0$, which is obtained from $\mathcal{H}_{\delta}$ by taking $\delta=0$. Specifically, $\mathcal{H}_0=(C,F_0,D,G_0)$, where $C$ and $D$ coincide with the definitions in \eqref{eqn:hybrid_system}, and the flow map $F_0$ and the jump map $G_0$ are given by
\begin{subequations}\label{eqn:nominal}
    \begin{align}
        \label{eqn:flow_map_nominal}
        &F_0(\xi) \coloneqq \{(A_qy + b_q)\} \times \{0\}\times\{0\} \times\{0\},\\
        \label{eqn:jump_map_nominal}
        &G_0(\xi) \coloneqq \{y\}\times(Q \backslash \{q\}) \times \{\tilde{u}\}\times\{\tau -1\}.
    \end{align}
\end{subequations}
Note that this system can exhibit at most one jump (i.e., one switch) because the timer $\tau$ is not allowed to increase, i.e., $\dot{\tau}=0$. Also, in this nominal system, the external load $\tilde{u}$ remains constant. {\color{black} Let} $R_0(K)$ denote the reachable set for the hybrid system $\mathcal{H}_0$ from the set of initial conditions $K$ (see Definition \ref{definitionReachable}), and consider the following set for each $q\in\mathcal{Q}$ and  $j\in\mathbb{Z}_{\geq 0}$:
\begin{align*}
\mathcal{A}_{j,q}=\overline{R_0 \left( \left ( \{y^*_q\} + (j+1)^{-1}\mathbb{B} \right ) \times \{q\} \times \mathcal{U} \times [0,1]\right )}.
\end{align*}
where, as before, $y^*_q:=-A_q^{-1}b_q$. Intuitively, the set \(\mathcal A_{j,q}\) characterizes the closure of the set of points that can be reached by trajectories of \(\mathcal H_0\) from initializations at the mode $q$, and with initial conditions of $y$ restricted to be in a ball of radius $(j+1)^{-1}$ around the $q^{th}$ equilibrium point. Next, define the sets
\begin{align}\label{eqn:omega_limit_set}
        \mathcal{A}_q  \coloneqq \bigcap_{j \in \mathbb{Z}_{\geq 0}} \mathcal{A}_{j,q},~~~~\text{and}~~~~~\mathcal{B} \coloneqq \bigcup _{q \in \mathcal{Q}} \mathcal{A}_q.
\end{align}
{\color{black} The above construction was introduced in \cite[Equation (21)]{baradaran2020omega}, without the additional set $\mathcal{U}$, where sequences of reachable sets were used to analyze switched differential inclusions subject to average dwell-time constraints (see Remark \ref{rem:prop_1}). Indeed, the proof of Proposition 4 in \cite{baradaran2020omega} can be adapted in a straightforward manner to show that the set $\mathcal{B}$ defined in \eqref{eqn:omega_limit_set} is the $\Omega$-limit set of the nominal hybrid system $\mathcal{H}_0$ from any sufficiently large compact set and, therefore, is the global attractor. However, the generality of the results in \cite{baradaran2020omega} precludes a simpler characterization of the $\Omega$-limit, i.e., the set $\mathcal{B}$. 
Herein, we aim to provide a sharper characterization of the set $\mathcal{B}$ by specializing to our setting.}
To that end, we have the following key Proposition.

\vspace{0.15cm}
\begin{prop}
        \label{prop:desc}
        { \color{black} The sets $\mathcal{B}$ and $\mathcal{A}$, defined in \eqref{eqn:omega_limit_set} and \eqref{eqn:omega_limit_set_1}, respectively, are equal.}
\end{prop}

\vspace{0.1cm}
\textbf{Proof:} To simplify notation, we begin by defining the following sets
\begin{align*}
        K^j_{q} &:=  (\{y^*_q\} + (j+1)^{-1} \mathbb{B}) \times \{q\} \times \mathcal{U} \times [0, 1], \\
        K^\infty_{q}&:= \{y^*_q\} \times \{q\} \times \mathcal{U} \times [0, 1].
\end{align*}
In this notation, we observe that
\begin{align*}
    \mathcal{A}_{j,q} = \overline{R_0 \big(K^j_{q}\big)}.
\end{align*}
We also let $\Theta^q_t:\mathbb{R}^{|\mathcal{G}|+2}\rightarrow \mathbb{R}^{|\mathcal{G}|+2}$ be the solution map of
\begin{align}\label{eq:individial_nominal_lti}
    \dot{y} = A_q y + b_q,
\end{align}
i.e. $\Theta^q_t(y_0)$ is the point corresponding to the unique solution of \eqref{eq:individial_nominal_lti} at time $t$ and starting from the initial condition $y(0)=y_0$. One can show that $\Theta^q_t$ coincides with the definition in \eqref{eq:flow_individual_linear_mode}. Since $K^\infty_{q}\subset K^j_{q}$, for all $j\geq 0$, it follows that 
    \begin{equation*}
        \overline{R_0 \left (K^\infty_{q}\right )} \subseteq \mathcal{A}_{{\color{blue}j,q}} ,
    \end{equation*}
    for all $j\in\mathbb{Z}_{\geq 0}$ and, therefore that 
    \begin{equation*}
        \overline{R_0 \left (K^\infty_{q}\right )} \subseteq \bigcap_{j\in\mathbb{Z}_{\geq 0}} \mathcal{A}_{{\color{blue}j,q}} = \mathcal{A}_{q}.
    \end{equation*}
    Next, we show the other direction of the inclusion. 

    As remarked above, the nominal hybrid system $\mathcal{H}_0$ can experience at most one jump. Therefore, for any $j\in\mathbb{Z}_{\geq 0}$, and for all maximal $\phi \in \mathcal{S}(K^{j}_{q})$, we have that one of the following cases holds:
    \begin{align*}
        (\text{C1})\quad\text{dom}(\phi) &= [0, \infty) \times \{0\}, \text{ or } \\
        (\text{C2})\quad\text{dom}(\phi) &= [0,t_1] \times \{0\} \cup [t_1,\infty) \times \{1\}, 
    \end{align*}
    for some $t_1 \geq 0$. In either case, we have that
    \begin{equation*}
        \phi(t,0) = (\Theta^q_t(y(0,0)), q, \tilde{u}(0,0), \tau(0,0)),
    \end{equation*}
    for all $(t,0) \in \text{dom}(\phi)$, where the initial conditions are such that $y(0,0) \in \left (\{y^*_q\} + (j+1)^{-1}\mathbb{B} \right )$, $\tilde{u}(0,0)\in\mathcal{U}$, and $\tau(0,0) \in [0,1]$. Due to Assumption \ref{asmp:hurwitz_Aq}, the equilibrium point $y^*_q$ is globally exponentially stable for the linear system \eqref{eq:individial_nominal_lti}. Hence, there exists a constant $c_1>0$ such that 
    $$\Theta^q_t(y(0,0)) \in \left (\{y^*_q\} + c_1(j+1)^{-1} \mathbb{B}  \right ),$$ 
    for all $(t,0) \in \text{dom}(\phi)$. Therefore, in the case (C1),
    \begin{align*}
        \phi(t,0)\in (\{y^*_q\} + c_1(j+1)^{-1} \mathbb{B}) \times \{q\}\times \mathcal{U} \times [0, 1],
    \end{align*}
    for all $(t,0)\in \text{dom}(\phi)$. 
    Suppose that the case (C2) holds. Then, for all $(t,0) \in \text{dom}(\phi)$,
    \begin{equation*}
        \phi(t,0) = (\Theta^q_t(y(0,0)), q, \tilde{u}(0,0), \tau(0,0)),
    \end{equation*}
    whereas, for all $(t,1) \in \text{dom}(\phi)$, we have that
    \begin{equation*}
        \phi(t,1) = (\Theta_{t-t_1}^{q^+}\circ\Theta_{t_1}^q(y(0,0)), q^+, \tilde{u}, {\color{blue} 0 } ),
    \end{equation*}
    for some $q^+\in\mathcal{Q}\backslash\{q\}$. Similar to the case (C1), there exists a constant $c_1>0$ such that
    \begin{equation*}
        \Theta^q_t(y(0,0)) \in \left (\{y^*_q\} + c_1(j+1)^{-1} \mathbb{B}  \right ),
    \end{equation*}
    for all $(t,0) \in \text{dom}(\phi)$. Therefore, we have that
    \begin{align*}
        \Theta_{t-t_1}^{q^+}\circ\Theta_{t_1}^q(y(0,0)) \in &\Big\{\Theta_{t-t_1}^{q^+}(y): \notag \\
        & ~~~~~~~|y-y^*_q|\leq c_1(j+1)^{-1}\Big\},
    \end{align*}
    for all $(t,1) \in \text{dom}(\phi)$. However, because the linear system \eqref{eq:individial_nominal_lti} is exponentially stable for all $q\in\mathcal{Q}$, for each $q\in\mathcal{Q}$, there exists a {\color{black} quadratic Lyapunov function  given by a positive definite matrix $P_q$. Consequently, there exists a constant $c_2>0$ such that
    \begin{align*}
        |\Theta_{t}^{q}(y)-\Theta_{t}^{q}(\tilde{y})| \leq c_2|y-\tilde{y}|,
    \end{align*}
    for all $y,\tilde{y}\in\mathbb{R}^{|\mathcal{G}|+2}$, all $t\geq0$, and all $q\in\mathcal{Q}$.} As a result, we arrive at the containment
    \begin{align*}
        \Big\{\Theta_{t-t_1}^{q^+}(y):~&|y-y^*_q|\leq c_1(j+1)^{-1}\Big\} \notag \\
        &\subseteq \Big\{\Theta_{t-t_1}^{q^+}(y^*_q)\Big\} + c_1 c_2 (j+1)^{-1}\mathbb{B}.
    \end{align*}
    {\color{blue}Since a jump requires $\tau=1$, and $\tau=0$ immediately after
jumps, we have
\begin{align}
R_0(K_q^\infty)
={}&
\{y_q^\ast\}\times\{q\}\times\mathcal{U}\times[0,1]
\nonumber\\
&\cup
\bigcup_{r\in\mathcal{Q}\setminus\{q\}}
\{\Theta_t^r(y_q^\ast):t\geq0\}
\times\{r\}\times\mathcal{U}\times\{0\}.
\label{eq:R0Kinf}
\end{align}
In general, the sets
$\{\Theta_t^r(y_q^\ast):t\geq0\}$ are not closed. Indeed,
since each subsystem is exponentially stable, we have $\lim_{t\to\infty}\Theta_t^r(y_q^\ast)=y_r^\ast$, and hence
\[
    \overline{\{\Theta_t^r(y_q^\ast):t\geq0\}}
    =
    \{\Theta_t^r(y_q^\ast):t\geq0\}
    \cup\{y_r^\ast\}.
\]
Consequently,
\begin{align}
&\overline{R_0(K_q^\infty)}
={}
\{y_q^\ast\}\times\{q\}\times\mathcal{U}\times[0,1]
\nonumber\\
&\cup
\bigcup_{r\in\mathcal{Q}\setminus\{q\}}
\left(
\{\Theta_t^r(y_q^\ast):t\geq0\}
\cup\{y_r^\ast\}
\right)
\times\{r\}\times\mathcal{U}\times\{0\}.
\label{eq:closureR0Kinf}
\end{align}
We now show that $\mathcal{A}_q=
    \overline{R_0(K_q^\infty)}$. Indeed, since $K_q^\infty\subset K_q^j$ for every
$j\in\mathbb{Z}_{\geq0}$, it follows that $R_0(K_q^\infty)
    \subseteq
    R_0(K_q^j)$, and therefore $\overline{R_0(K_q^\infty)}
    \subseteq
    \overline{R_0(K_q^j)}
    =
    \mathcal{A}_{j,q}$ for every $j\in\mathbb{Z}_{\geq0}$. Thus,
\begin{equation}
    \overline{R_0(K_q^\infty)}
    \subseteq
    \bigcap_{j\in\mathbb{Z}_{\geq0}}\mathcal{A}_{j,q}
    =
    \mathcal{A}_q.
\label{eq:first-inc}
\end{equation}
To establish the reverse inclusion, define $c_3:=\max\{c_1,c_1c_2\}$. Comparing each solution from $K_q^j$ with the solution from
$K_q^\infty$ having the same initial mode, load, timer, jump
time, and target mode, we obtain
\[
    R_0(K_q^j)
    \subseteq
    R_0(K_q^\infty)
    +c_3(j+1)^{-1}\mathbb{B},
\]
Taking closures gives
\[
    \mathcal{A}_{j,q}
    =
    \overline{R_0(K_q^j)}
    \subseteq
    \overline{
    R_0(K_q^\infty)
    +c_3(j+1)^{-1}\mathbb{B}}.
\]
Since $R_0(K_q^\infty)
    \subseteq
    \overline{R_0(K_q^\infty)}$ we further obtain $\mathcal{A}_{j,q}
    \subseteq
    \overline{R_0(K_q^\infty)}
    +c_3(j+1)^{-1}\mathbb{B}$, where we have used the fact that
$\overline{R_0(K_q^\infty)}$ is closed and the ball is compact.
Therefore,
\begin{align}
    \mathcal{A}_q
    &=
    \bigcap_{j\in\mathbb{Z}_{\geq0}}\mathcal{A}_{j,q}
    \nonumber\\
    &\subseteq
    \bigcap_{j\in\mathbb{Z}_{\geq0}}
    \left(
    \overline{R_0(K_q^\infty)}
    +c_3(j+1)^{-1}\mathbb{B}
    \right)
    \nonumber\\
    &=
    \overline{R_0(K_q^\infty)}.
\label{eq:second-inc}
\end{align}
Combining \eqref{eq:first-inc} and \eqref{eq:second-inc} yields
$\mathcal{A}_q=\overline{R_0(K_q^\infty)}$. Finally, taking the union over $q\in\mathcal{Q}$ and using
\eqref{eq:closureR0Kinf}, we obtain
\begin{align*}
&\mathcal{B}
={}
\bigcup_{q\in\mathcal{Q}}
\overline{R_0(K_q^\infty)}
\\
&={}
\bigcup_{q\in\mathcal{Q}}
\left[
\{y_q^\ast\}\times\{q\}\times\mathcal{U}\times[0,1]
\right.\\[-1mm]
&\left.
{}\cup
\bigcup_{r\in\mathcal{Q}\setminus\{q\}}
\left(
\{\Theta_t^r(y_q^\ast):t\geq0\}
\cup\{y_r^\ast\}
\right)
\times\{r\}\times\mathcal{U}\times\{0\}
\right].
\end{align*}
For every $r\in\mathcal{Q}$, the additional limiting set $\{y_r^\ast\}\times\{r\}\times\mathcal{U}\times\{0\}$ is already contained in $\{y_r^\ast\}\times\{r\}\times\mathcal{U}\times[0,1]$. Hence,
\begin{align*}
\mathcal{B}
={}&
\bigcup_{q\in\mathcal{Q}}
\left[
\{y_q^\ast\}\times\{q\}\times\mathcal{U}\times[0,1]
\right.\\[-1mm]
&\qquad\left.
{}\cup
\bigcup_{r\in\mathcal{Q}\setminus\{q\}}
\{\Theta_t^r(y_q^\ast):t\geq0\}
\times\{r\}\times\mathcal{U}\times\{0\}
\right]
\\
={}&\mathcal{A}.
\end{align*}
This completes the proof.}
    \hfill $\blacksquare$

\vspace{0.1cm}
\begin{rem}
{\color{black}The result of Proposition \ref{prop:desc} provides a sharp characterization of the $\Omega$-limit set for the hybrid system of interest. Similar characterizations are presented in \cite{baradaran2020omega} for more general nonlinear switched systems that can have multiple consecutive jumps. In this sense, Proposition 1 leverages the additional affine structure of \eqref{eqn:nominal}, as well as the fact that the system can switch at most once before flowing again, to derive a sharper expression for the resulting $\Omega$-limit set. To the best of our knowledge, this expression has not been explicitly  presented before in the literature.}  
\end{rem}
With Proposition \ref{prop:desc} 
, we can now prove Theorem \ref{thm:stability_result}.

\vspace{0.1cm}
\noindent 
\textbf{Proof of Theorem \ref{thm:stability_result}:} The proof follows similar steps as in \cite[Sec. 7]{goebel2012hybrid} and \cite{baradaran2020omega}, but leveraging the set $\mathcal{A}$ characterized in Proposition \ref{prop:desc}. In particular, by construction the hybrid system $\mathcal{H}_0$ satisfies the hybrid basic conditions. Let $\Omega(K)$ be the $\Omega$-limit set of the hybrid system $\mathcal{H}_0$ from the set $K$ where $K$ is chosen to be sufficiently large such that the set $\mathcal{A}$ is contained in the interior of $K$. It follows from \cite[Prop. 4]{baradaran2020omega} and Proposition \ref{prop:desc} that $\Omega(K)=\mathcal{A}$, and by \cite[Corollary 7.7]{goebel2012hybrid} we obtain that $\Omega(K)$ is asymptotically stable with basin of attraction containing $K$. Since $K$ can be taken arbitrarily large, the set $\Omega(K)$ is globally asymptotically stable for the hybrid system $\mathcal{H}_0$. 

Next, consider the inflated hybrid system $\mathcal{H}_{\epsilon}$
where $\epsilon>0$, and where the data of \eqref{eqn:inflation} is constructed from the data of $\mathcal{H}_0$ {\color{black} as in Definition \ref{defn:inflation}}. By \cite[Thm. 7.21]{goebel2012hybrid}, the inflated hybrid system $\mathcal{H}_{\epsilon}$ renders the set $\mathcal{A}$ semi-globally practically asymptotically stable as $\epsilon\to0^+$. Moreover, since $\Pi(\cdot)$ is OSC and LB, and $\mathcal{U}$ and $\mathcal{Q}$ are compact, it follows that there exists $L>0$ such that $$\sup_{\substack{
q \in \mathcal{Q} \\
\tilde{u} \in \mathcal{U} \\
v \in \Pi(\tilde{u})
}}
\left\{
|v|,\,
\left|A_q^{-1}B_q v\right|
\right\}
\leq L$$ for all $\tilde{u}\in\mathcal{U}$ and all $q\in\mathcal{Q}$. For each $\epsilon$, let $\delta$ be sufficiently small such that $\max\{\delta_1, \delta_2 L, \delta_3\} < \epsilon$. It follows that for each $\epsilon>0$ the data of the hybrid system \eqref{eqn:hybrid_system} is contained on the data of $\mathcal{H}_{\epsilon}$, and every solution of \eqref{eqn:hybrid_system} is also a solution to \eqref{eqn:inflation}. The latter fact implies that system \eqref{eqn:hybrid_system} also renders $\mathcal{A}$ semi-globally practically asymptotically stable as $|\delta|\to0^+$.  

\vspace{-0.2cm}
\subsection{Arbitrarily Varying Loads}\label{subsec:approach_2}
In the previous section, we established practical stability for slowly varying loads. However, real-world loads may exhibit arbitrarily fast variations over bounded time intervals. In this section, we study such scenario via input-to-state stability tools for hybrid systems, and we analyze stability with respect to a larger set than that in Theorem 1. 

In a manner similar to {\color{blue} S}ubsection \ref{subsec:approach_1}, we model the switching behavior via the hybrid automaton \eqref{eq:swtching_automata}. We also introduce the following change of coordinates
\begin{align}
    y = x+ A_q^{-1}b_q.
\end{align}
Direct differentiation shows that the continuous time evolution of the state $y$ is given by $\dot{y}= A_q y + B_q \tilde{u}$. On the other hand, whenever the mode $q$ switches to $q^+\in\mathcal{Q}\backslash\{q\}$, the jumps of $y$ satisfy $y^+= y +A_{q^+}^{-1}b_{q^+}- A_{q}^{-1}b_{q} $. Therefore, the combined state $\xi=(y,q,\tau)$ evolves according to the hybrid system with input $\mathcal{H}_{\tilde{u}}=(C,\tilde{F},D,\tilde{G})$ \cite{cai2009characterizations}, defined by the flow and jump sets
\begin{subequations}\label{eqn:modified_flow_jump}
\begin{align*}
        \!\!\!\!\!\!\!\!C\coloneqq \mathbb{R}^{\lvert \mathcal{G} \rvert + 2} \times \mathcal{Q} \times [0,N_0],~D \coloneqq \mathbb{R}^{\lvert \mathcal{G} \rvert + 2} \times \mathcal{Q} \times [1,N_0],
\end{align*}
where $N_0\geq 1$ is the chatter bound (i.e. the maximum allowable number of consecutive jumps), and the flow and jump maps
\begin{align}
    \tilde{F}(\xi,\tilde{u})&:= \{A_q y + B_q\tilde{u}\}\times\{0\}\times[0,\eta], \\
    \tilde{G}(\xi,\tilde{u})&:=\bigcup_{r\in\mathcal{Q}\backslash\{q\}}\{y +A_{r}^{-1}b_{r}- A_{q}^{-1}b_{q} \}\times\{r\} \notag\\
    &~~~~~~~~~~~~~~~~\times\{\tau-1\},
\end{align}
\end{subequations}
where $\eta\geq 0$ controls the rate of change of the state of the timer $\tau$. As shown in \cite[Ex. 2.15]{goebel2012hybrid}, for every $0\leq s<t$, every solution of this system satisfies the ADT bound \eqref{eq:ADT} with $\tau_d=1/\eta$.
Moreover, every switching signal satisfying the ADT bound \eqref{eq:ADT} can be generated using \eqref{eqn:modified_flow_jump} by selecting the appropriate initial condition. 
We now present the second main result of the paper which concerns the input-to-state stability properties of the hybrid system $\mathcal{H}_{\tilde{u}}$ with respect to the input $\tilde{u}$.
\begin{thm}\label{thm:ISS_result}
    Let Assumption \ref{asmp:hurwitz_Aq} be satisfied, and define the set $\mathcal{A}_0:=\{0\}\times\mathcal{Q}\times[0,N_0]$. Then, there exists $\eta^*\in(0,\infty)$ such that, for any $\eta\in(0,\eta^*)$, there exists positive constants $\kappa_i>0$, $i\in\{1,2,3\}$, such that, for any continuous input $\tilde{u}:[0,\infty)\rightarrow \mathbb{R}^{|\mathcal{G}|+1}$ with $|\tilde{u}|_\infty < \infty$, every maximal solution of the hybrid system defined by \eqref{eqn:modified_flow_jump} is complete and satisfies the uniform bound
    \begin{align}\label{eq:ISS_exp_KL}
        \!\!\!\!\!\!\!\!|\xi(t,j)|_{\mathcal{A}_0}\leq \kappa_1 \mathrm{e}^{-\kappa_2(t+j)}|\xi(0,0)|_{\mathcal{A}_0} + \kappa_3( |\tilde{u}|_{\infty} + c),
    \end{align}
    for all $(t,j)\in \mathrm{dom}(\xi)$, where $c\geq0$ is defined by
    \begin{align}
        c:= \max_{q\in\mathcal{Q}} \max_{r\in\mathcal{Q}\backslash\{q\}} |A_{r}^{-1}b_{r}-A_q^{-1}b_q|.
    \end{align}
\end{thm}

\vspace{0.1cm}
\begin{rem} 
Theorem \ref{thm:ISS_result} establishes the existence of a lower bound on the dwell-time of the switching signal such that, for any larger dwell-time, the frequency deviation, mechanical powers, and the state of secondary controller, evolving according to the hybrid system defined by \eqref{eqn:modified_flow_jump} settle in a certain neighborhood of the set $\mathcal{A}$. Moreover, the theorem establishes an upper bound on the size of this neighborhood that is directly proportional to the magnitude of the time-varying load $\tilde{u}$. In particular, the larger the magnitude of the time-varying load, the more conservative the estimate on the convergence of the trajectories becomes. In the absence of the exogenous input $\tilde{u}$, the ultimate
bound is determined by the nonnegative constant $c$, which measures the maximum distance between the mode-dependent equilibria. In particular, this bound vanishes when all mode-dependent equilibria coincide. 
\hfill $\Box$
\end{rem}

\noindent 
\textbf{Proof of Theorem 2:}
Let $\hat{u}:=(u_c,u_d)$, and consider the hybrid system with inputs, given by
\begin{equation}\label{eqn:hybrid_system_input}
   \mathcal{H}_{\hat{u}}:\begin{cases}
        (\xi,\hat{u}) \in C \times \hat{\mathcal{U}}, & \dot{\xi}\in F(\xi,\hat{u})\\ 
        (\xi,\hat{u}) \in D \times \hat{\mathcal{U}}, & \xi^+\in G(\xi,\hat{u}),
    \end{cases}
\end{equation}
where the set of admissible inputs $\hat{\mathcal{U}}$ is given by
\begin{align}
    \hat{\mathcal{U}}:=\mathbb{R}^{(1+|\cG|)}\times\bigcup_{q\in\mathcal{Q}}\bigcup_{r\in\mathcal{Q}\backslash\{q\}}\{A_{r}^{-1}b_{r}-A_{q}^{-1}b_{q}\},
\end{align}
the flow and jump sets coincide with \eqref{eqn:modified_flow_jump}, and the flow and jump maps are
\begin{align}
    F(\xi,\hat{u})&:=\{A_q y + B_q u_c\}\times\{0\}\times[0,\eta],\\
    G(\xi,\hat{u})&:=\{y + u_d\}\times(\mathcal{Q}\backslash\{q\})\times\{\tau-1\}.
\end{align}
Then, every trajectory of the hybrid system $\mathcal{H}_{\tilde{u}}$ coincides with a trajectory of the hybrid system $\mathcal{H}_{\hat{u}}$ for an appropriate choice of the hybrid input $\hat{u}$. Indeed, if $\xi$ is a solution of the hybrid system $\mathcal{H}_{\tilde{u}}$ with a corresponding input $\tilde{u}$, then the hybrid arc $\tilde{\xi}$ with input $\hat{u}$ defined by
\begin{align*}
\tilde{\xi}(t,j) &:= \xi(t,j),~~~~u_c(t,j) := \tilde{u}(t),\\
    u_d(t,j) &:= A_{q(t,j+1)}^{-1}b_{q(t,j+1)}-A_{q(t,j)}^{-1}b_{q(t,j)} ,
\end{align*}
is a solution to the hybrid system $\mathcal{H}_{\hat{u}}$. Note that there is no violation of causality in this choice of the input since the input $\hat{u}$ is specified as a hybrid signal, not as a function of the state of the hybrid system $\mathcal{H}_{\hat{u}}$. From the preceding discussion, we see that stability properties of $\mathcal H_{\hat u}$ with respect to
$\hat u$ entail corresponding practical stability properties of
$\mathcal H_{\tilde u}$ with respect to $\tilde u$. Therefore, we focus our attention on the hybrid system $\mathcal{H}_{\hat{u}}$. 

By Assumption \ref{asmp:hurwitz_Aq}, there exist symmetric positive definite matrices $P_q\succ 0$ such that, for all $q\in\mathcal{Q}$, we have that
\begin{align}
	Q_q:= -(A_q^\top P_q + P_q A_q) \succ 0.
\end{align}
\begin{figure*}[t!]
\centering\includegraphics[width=0.33\linewidth]{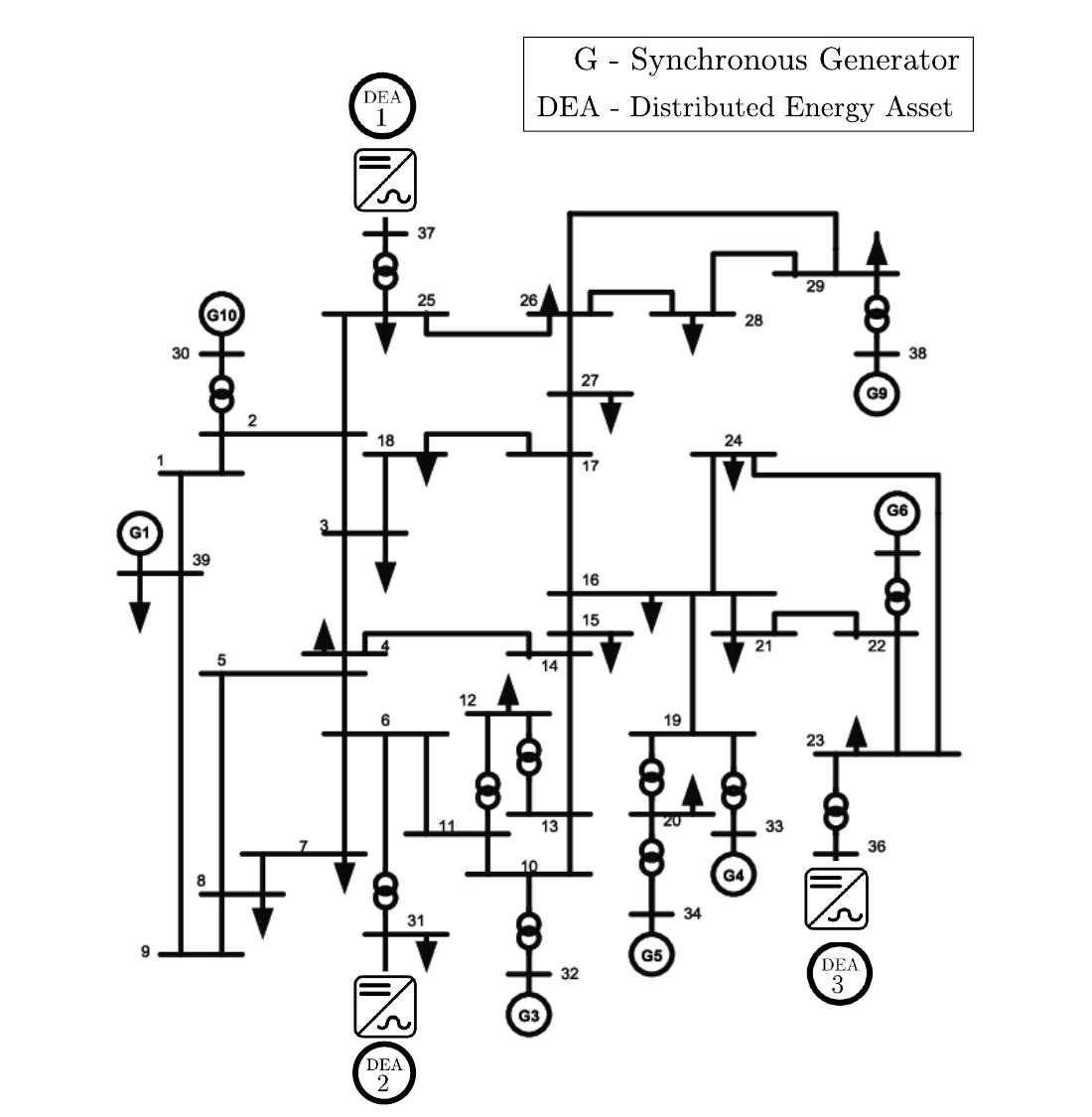}\includegraphics[width=0.33\linewidth]{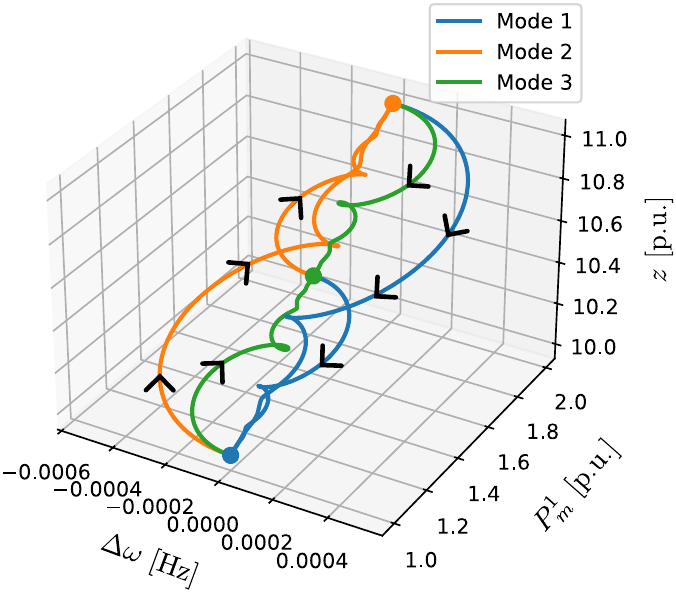}\includegraphics[width=0.33\linewidth]{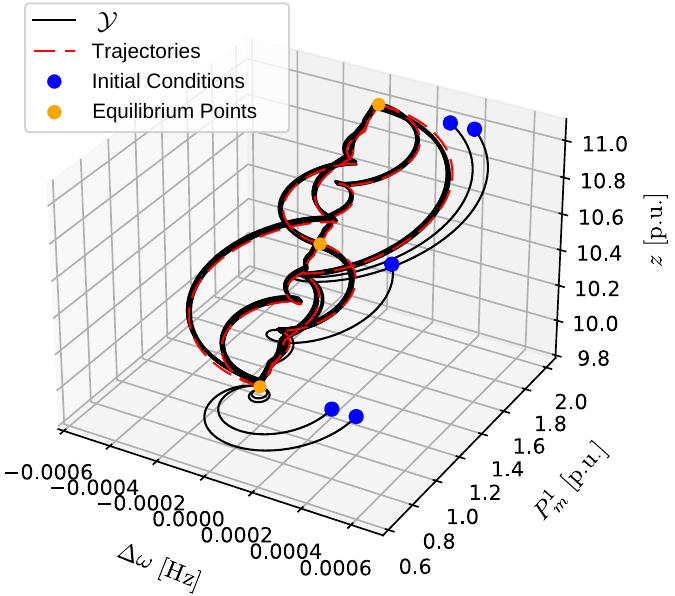}
    \caption{Left: Line diagram of the IEEE-39 bus system with synchronous generators and DEAs. Center: Projection onto the $(\Delta \omega, P^1_m, z)$ axis of the set $\mathcal{Y}$ for the IEEE 39-bus system . Right: Trajectories of the nominal 39-bus system under sufficiently slowly switching loads.}
    \label{fig:ieee_line}
\end{figure*}
For any $\mu \geq 0$ and any choice of the matrices $P_q$, we define the function $V:C\cup D\rightarrow\mathbb{R}_{\geq 0}$ by
\begin{align}
	V(\xi)= \mathrm{e}^{\mu \tau} y^\top P_q y,
\end{align}
which satisfies
\begin{align}\label{eq:PD_radial_unbounded}
\underline{\alpha}|\xi|_{\mathcal{A}_0}^2\leq V(\xi)\leq \overline{\alpha}|\xi|_{\mathcal{A}_0}^2, \quad \forall \xi\in C\cup D\cup G(D)
\end{align}
where $\underline{\alpha}:=\min_{q\in\mathcal{Q}}\sigma_{\min}(P_q)$ and $ \overline{\alpha}:=\mathrm{e}^{\mu N_0} \max_{q\in\mathcal{Q}}\sigma_{\max}(P_q)$. Direct computation shows that, for all $(\xi,\hat{u})\in C\times\hat{\mathcal{U}}$, the function $V$ satisfies the inequality
\begin{align}\label{eq:ISS_V_flow_1}
	\!\!\!\!\!\!\!\!\dot{V}(\xi,\hat{u}) \leq -\mathrm{e}^{\mu\tau} y^\top (Q_q - \eta \mu P_q ) y + 2 \mathrm{e}^{\mu\tau} y^\top P_q B_q u_c,
\end{align}
where $\dot{V}(\xi,\hat{u}) = \max_{f \in F(\xi,\hat{u})}\langle\nabla V,f\rangle$. By adding and subtracting terms, the inequality \eqref{eq:ISS_V_flow_1} can be rewritten as
\begin{align}
	\dot{V}(\xi,\hat{u}) \leq  &-\mathrm{e}^{\mu\tau}[y^\top, u_c^\top] S [y^\top, u_c^\top]^\top \notag\\
    &- (1-\theta) \mathrm{e}^{\mu\tau} y^\top Q_q y + \mathrm{e}^{\mu\tau} u_c^\top R u_c,
\end{align}
for any $\theta\in(0,1)$ and any symmetric positive definite matrix $R$, where the matrix $S$ is given by
\begin{align}
	S = \begin{bmatrix}
		\theta Q_q - \eta \mu P_q & -P_q B_q \\
		- B_q^\top P_q & R
	\end{bmatrix}.
\end{align}
Since $Q_q\succ 0$ for all $q\in\mathcal{Q}$, it follows that, for any $\theta\in(0,1)$ and any $\mu \geq 0$, there exists $\eta^*>0$ such that 
\begin{align}
    \theta Q_q - \eta \mu P_q \succ 0,
\end{align}
for all $\eta\in(0,\eta^*)$ and all $q\in\mathcal{Q}$. Moreover, for any $\eta\in(0,\eta^*)$, there exists $R\succ 0$ such that $S\succeq 0$, for all $q\in\mathcal{Q}$. Therefore, for any such choice of $\eta$ and $R$, the function $V$ satisfies the inequality
\begin{align}
	\dot{V}(\xi,\hat{u}) \leq - \mathrm{e}^{\mu\tau} (1-\theta) y^\top Q_q y + \mathrm{e}^{\mu N_0} u_c^\top R u_c.
\end{align}
On the other hand, for all $(\xi,\hat{u})\in D\times\hat{\mathcal{U}}$
and every $\xi^+\in G(\xi,\hat{u})$, we have
$y^+=y+u_d$ and $\tau^+=\tau-1$. Hence,
\begin{align}
V(\xi^+)-V(\xi)
={}&
e^{\mu\tau}\Big[
y^\top\big(e^{-\mu}P_{q^+}-P_q\big)y
+2e^{-\mu}y^\top P_{q^+}u_d
\notag\\
&\hspace{32mm}
+e^{-\mu}u_d^\top P_{q^+}u_d
\Big]
\notag\\
={}&
-e^{\mu\tau}
\begin{bmatrix}
y^\top & u_d^\top
\end{bmatrix}
Y
\begin{bmatrix}
y^\top & u_d^\top
\end{bmatrix}^{\!\top}
\notag\\
&-(1-\theta)e^{\mu\tau}y^\top P_qy
+e^{\mu\tau}u_d^\top Mu_d,
\end{align}
where{\color{blue}
\[
Y:=
\begin{bmatrix}
\theta P_q-e^{-\mu}P_{q^+}
&
-e^{-\mu}P_{q^+}
\\
-e^{-\mu}P_{q^+}
&
M-e^{-\mu}P_{q^+}
\end{bmatrix}.
\]
Let $\mu>0$ be sufficiently large, such that
\[
H_{q,q^+}:=\theta P_q-e^{-\mu}P_{q^+}\succ0
\]
for every $q,q^+\in\mathcal{Q}$. Finally, let $M\succ0$ be
sufficiently large, such that
\[
M-e^{-\mu}P_{q^+}
-e^{-2\mu}P_{q^+}H_{q,q^+}^{-1}P_{q^+}
\succeq0
\]
for every $q,q^+\in\mathcal{Q}$.} By the Schur complement,
this guarantees that $Y\succeq0$.
%
Hence, for any such choices, the function $V$ satisfies the inequality
\begin{align}
	\Delta V(\xi,\hat{u})&\leq - \mathrm{e}^{\mu \tau}(1-\theta) y^\top P_q y + \mathrm{e}^{\mu  N_0}u_d^\top M u_d,
\end{align}
where $\Delta V(\xi,\hat{u}):=\max_{\xi^+\in G(\xi,\hat{u})} V(\xi^+)-V(\xi)$. 

Defining the constants
\begin{align}
	\lambda&:= (1-\theta)\min\left\{\min_{q\in\mathcal{Q}}\frac{\sigma_{\min}(Q_q)}{\sigma_{\max}(P_q)},1\right\},\\ 
	\rho&:= \mathrm{e}^{\mu  N_0} \max \{\sigma_{\max}(R), \sigma_{\max}(M)\},
\end{align}
and after some algebraic manipulations, we arrive at
\begin{align}
        \label{eq:ISS_decrease_flow}
	\!\!\!\!\!\!\!\!\dot{V}(\xi,\hat{u}) &\leq - \lambda V(\xi)+ \rho |u_c|^2, \quad \forall (\xi,\hat{u})\in C\times \hat{\mathcal{U}}, \\
        \label{eq:ISS_decrease_jump}
	\!\!\!\!\!\!\!\!\Delta V(\xi,\hat{u}) &\leq - \lambda V (\xi) + \rho |u_d|^2, \quad \forall (\xi,\hat{u})\in D\times \hat{\mathcal{U}}.
\end{align}
Combining \eqref{eq:PD_radial_unbounded} and \eqref{eq:ISS_decrease_flow}-\eqref{eq:ISS_decrease_jump}, and using 
\cite[Lemma 9]{ochoa2024prescribed}, there exists positive constants $\kappa_i$, $i\in\{1,2,3\}$, such that every trajectory of the hybrid system $\mathcal{H}_{\hat{u}}$ satisfies the uniform exponential bound
\begin{align}
    |\xi(t,j)|_{\mathcal{A}_0}\leq \kappa_1 \mathrm{e}^{-\kappa_2(t+j)}|\xi(0,0)|_{\mathcal{A}_0} + \kappa_3 |\hat{u}|_{(t,j)},
\end{align}
where $|\hat{u}|_{(t,j)}$ is a shorthand notation for
\begin{align*}
    |\hat{u}|_{(t,j)}
:= \sup_{\substack{(0,0)\preceq(s,i)\preceq(t,j)\\
(s,i)\in\operatorname{dom}(\hat{u})}}
|\hat{u}(s,i)|.
\end{align*}
In addition, every maximal solution of the hybrid system is complete. On the other hand, we observe that any solution of the hybrid system $\mathcal{H}_{\hat{u}}$ is such that $\hat{u}(t,j)\in\hat{\mathcal{U}}$. Therefore, using the triangle inequality, we arrive at the upper bound
\begin{align*}
    |\hat{u}(t,j)| \leq |u_c(t)| + \underbrace{\max_{q\in\mathcal{Q}} \max_{r\in\mathcal{Q}\backslash\{q\}} |A_{r}^{-1}b_{r}-A_q^{-1}b_q|}_{=c}.
\end{align*}
We also recall that every solution of the hybrid system $\mathcal{H}_{\tilde{u}}$ defined in \eqref{eqn:modified_flow_jump} coincides with a solution of the hybrid system $\mathcal{H}_{\hat{u}}$ for some choice of the input $\hat{u}$. Moreover, the same conclusion regarding completeness of solutions can be drawn with respect to the hybrid system $\mathcal{H}_{\tilde{u}}$ defined in \eqref{eqn:modified_flow_jump}. Combining all of the above, we arrive at the bound in \eqref{eq:ISS_exp_KL} which concludes the proof. \hfill $\blacksquare$
%
\section{Numerical Experiments in the IEEE 39-Bus System}
\label{sec:sims}

\noindent
We  consider the IEEE 39-bus test system, which features 10 CGs. As shown in the left plot of Figure \ref{fig:ieee_line}, we added three large-scale DEAs. The other parameters for the generators are taken from the IEEE 39-bus test system data. 

\subsection{Aggregate model} We start with numerical experiments with the aggregate model~\eqref{eqn:switched_loads}.  
For the CGs, we set $\tau_g = 2, D_g = 1.5$ and $K_{\mathrm{gov},g} = \frac{1}{0.05},\ \forall g \in \cG$. For the DEAs, we set $M_d \in \{40,30,25\} , D_d \in \{2,1,3\},\ \forall d \in \cD$. For the secondary controller,  $\tau_z = 10$, $\beta = -0.1$ and $\zeta = \{\zeta_i\}_{i\in\mathcal{G}}$ such that every generator participates equivalently in the secondary frequency response. By leveraging Theorem 1 and  Proposition \ref{prop:desc} we compute the $\Omega$-limit set by running simulations of the three subsystems from given initial conditions. 
The center plot in Figure \ref{fig:ieee_line} shows the $y$-component of the $\Omega$-limit set; more specifically, it is the projection onto the axes of the frequency deviation $(\Delta \omega)$, first mechanical power $(P_m^1)$, and the secondary controller $(z)$. It was plotted using the explicit construction given in Proposition \ref{prop:desc}. The effect of the flows of the individual subsystems can clearly be seen. For instance, in orange, we observe two trajectories converging exponentially fast to the unique equilibrium point of subsystem $2$, starting from the equilibrium points of subsystems $1$ (blue) and $3$ (green). Note that for using the definition given in Equation \eqref{eqn:omega_limit_set} to plot the limit-set for the power transmission system may have been cumbersome. The right plot in Figure \ref{fig:ieee_line} illustrates the state trajectories of the IEEE 39-bus system under slowly switching loads. In color red (dashed), we plot again the  projection of the $\Omega$-limit set, while in color black (solid), we show the state trajectories of the actual system initialized from $5$ random initial conditions onto the $(\Delta \omega, P^1_m, z)$ axis. The actual system trajectories of the power system closely approximate the set $\mathcal{Y}$ --- a projection onto the $y$-component of the $\Omega$-limit set. We would like to emphasize that the actual system incorporates slow switching, and slowly-varying additive disturbances to the flow and jump dynamics of the nominal system. Figure \ref{fig:actual_39_time} shows the trajectories of the hybrid system \eqref{eqn:hybrid_system} for non-zero but sufficiently small value of $\delta$. The frequency deviation is driven to zero with minute deviations corresponding exactly to the switching instants. The switching signal causing this behavior can be seen in the bottom right of the same figure. To drive the frequency deviations to zero, the secondary controller oscillates in the interval $[9.5,11.5]$. Lastly, the mechanical power $P^m_1$ also oscillates in response to the oscillating reference signal $P^r_{\mathcal{G}}$. The oscillations of the frequency deviation, secondary controller, and the first mechanical power are best visualized in phase space. Lastly, Figure \ref{fig:distance} shows the distances between the system trajectories of the actual system and the set $\mathcal{Y}$ of the nominal system for all five initial conditions seen in the right plot of Figure \ref{fig:ieee_line}. The distances to the set of interest converge to a small neighborhood of zero, further validating the SGPAS result proved in Theorem \ref{thm:stability_result} for the state trajectories of the power system. 

{\color{black} The simulations so far show that the system trajectories of IEEE 39-bus system closely approximate the projection of the Omega-limit set ({the set $\mathcal{Y}$) of the ideal system $\mathcal{H}_0$. Theorem \ref{thm:ISS_result} however, claims that the system trajectories converge to only a neighborhood of the set $\mathcal{A}$ whose magnitude depends on both the input size and the maximum distance between the equilibrium points. While we acknowledge this conservatiness, we would like to mention that it is a direct result of the presence of multiple equilibria and the choice of a Lyapunov-based analysis. This can be seen clearly from the proof of Theorem \ref{thm:ISS_result} where the change of variables results in non-trivial jump dynamics in System \ref{eqn:hybrid_system_input}. This persistent deviation at every jump is precisely what gives rise to the constant $c$ seen in \eqref{eq:ISS_exp_KL}. }

\begin{figure}
    \centering
\includegraphics[width=\linewidth]{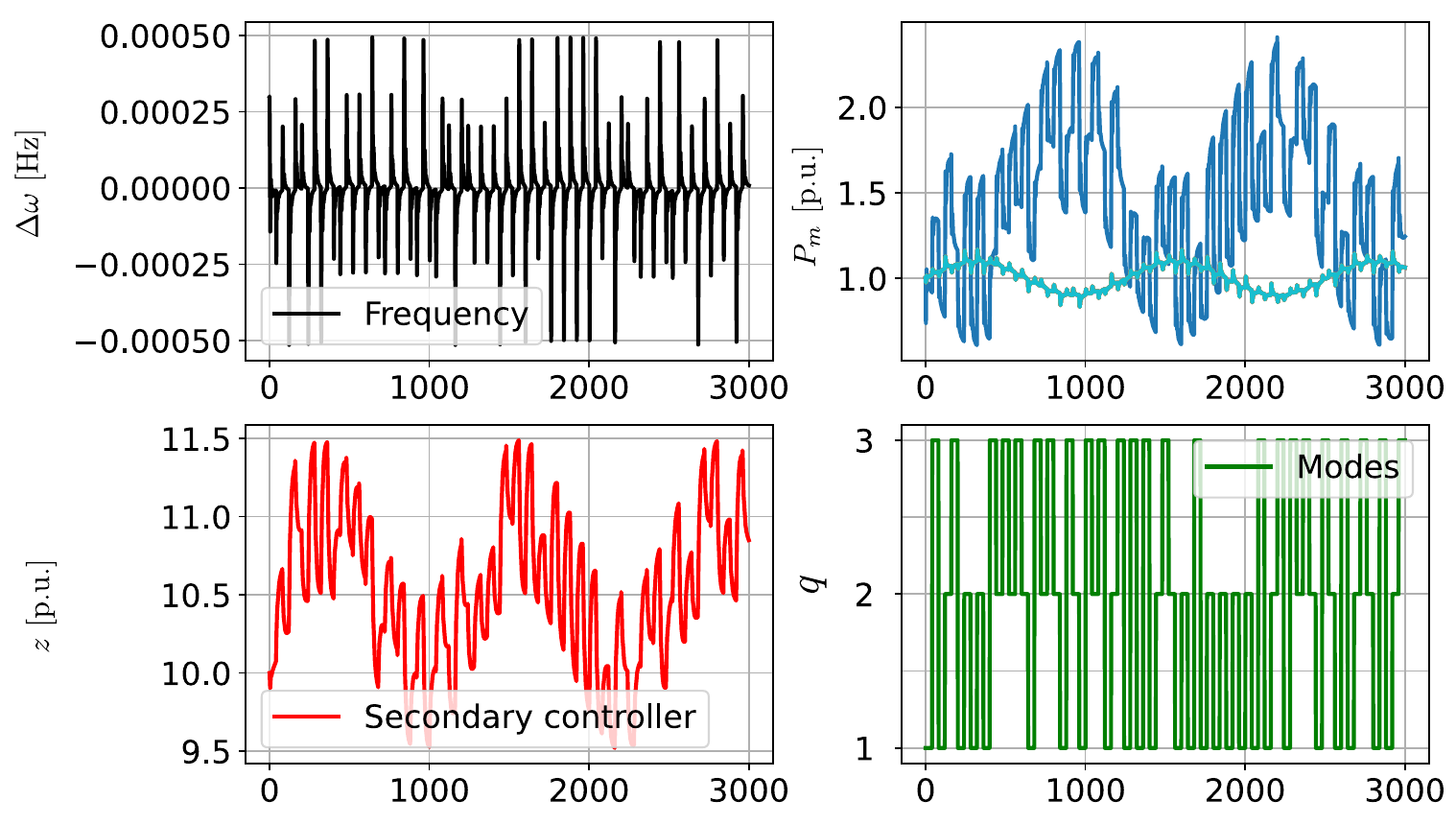}
\vspace{-0.5cm}
    \caption{Trajectories of the IEEE 39-bus system under switching loads and slowly varying sinusoidal loads.}
    \label{fig:actual_39_time}
\end{figure}
\begin{figure}[t!]
    \centering
    \includegraphics[width=0.9\linewidth]{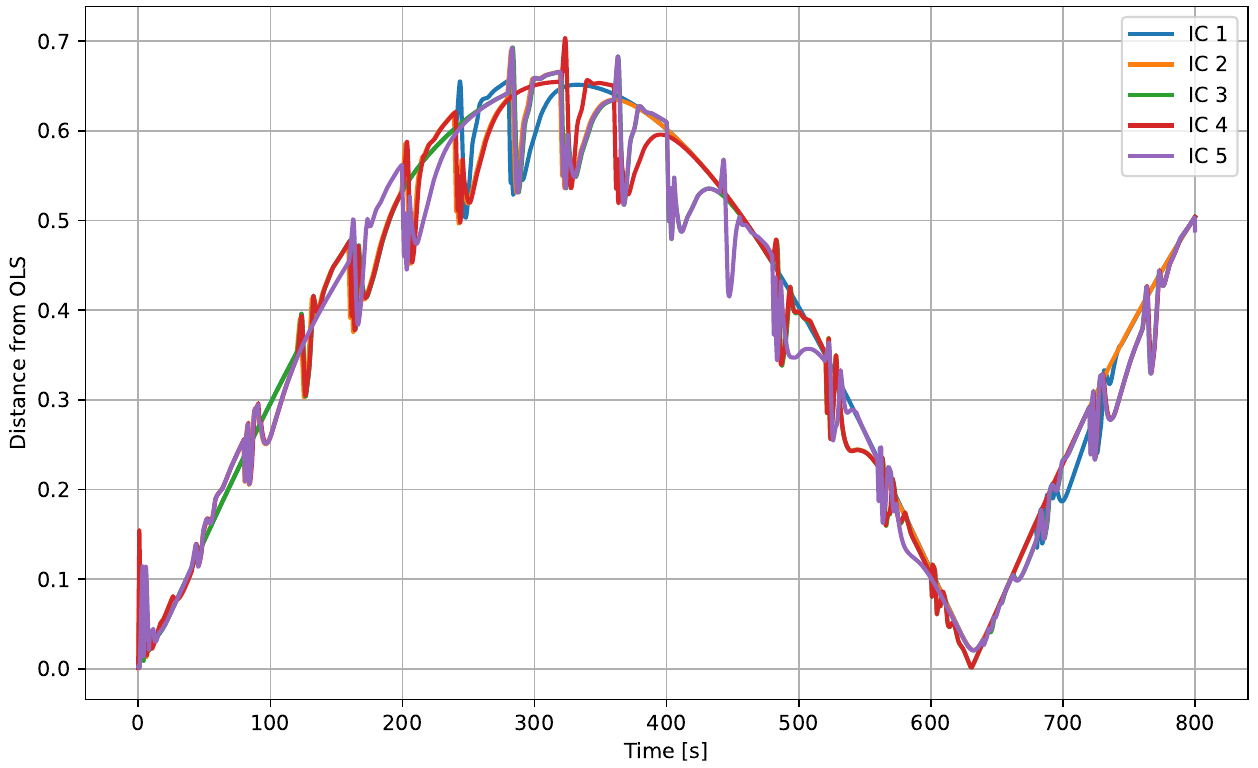}
    \caption{Distance of  power system state trajectories from the {\color{blue} set $\mathcal{Y}$} corresponding to the IEEE 39-bus system subject to sinusoidal loads.}
    \label{fig:distance}
    \vspace{-0.5cm}
\end{figure}
%

\vspace{0.2cm}
\noindent

\subsection{Model in~\eqref{eqn:ss_full}} To further validate our results and show the applicability of our findings to a  system with higher fidelity, we now consider the model in~\eqref{eqn:ss_full}. We add an integral action $K_i \delta_g$ in the right hand side of~\eqref{eq:gen_p} to steer the frequency deviation to zero even without the secondary controller. Unfortunately, the matrix $\bar{A}$ in~\eqref{eq:gen_p} has one zero eigenvalue; this is consistent 
with the rotational invariance of the angles in the power flow equations. We perform a change of coordinates $x = T x_\text{r}$, where $T$ is built based on the eigenvectors associated with the non-zero eigenvalues (see, e.g.,~\cite{colombino2019online}) and obtain a system $\dot x = Ax + Bu$ where $A = T^\top \bar{A} T$ and $B = T^\top \bar{B}$. Here, the matrix $A$ is Hurwitz. 

We consider three modes corresponding to three different inertia matrices $M_1, M_2, M_3$, damping matrices $D_1, D_2, D_3$, and loads $P_{\text{load}}^1, P_{\text{load}}^2, P_{\text{load}}^3$. The setpoints for CGs $P_{\mathcal{G}}^*$ remains unchanged. In Figure \ref{fig:nominal_39_full}, we see the projection of the $\Omega$-limit set of the full-order model of the IEEE 39-bus system onto the frequency deviation ($\Delta \omega $), the first ($P^m_1$), and the second ($P^m_2$) mechanical powers. We have chosen three modes for simplicity. As suggested by Proposition \ref{prop:desc}, we see the trajectories of nominal system under the flow of the other subsystems towards the corresponding equilibrium point. For instance, in blue, we see two trajectories converging exponentially fast to the corresponding equilibrium point (also in blue). Thus, even for the full-order model, with considerably more states than the reduced-order model ($65$ versus $12$), Proposition \eqref{prop:desc} gives a computationally tractable scheme for visualizing the asymptotic behavior of the IEEE 39-bus system. 

Figure \ref{fig:actual_39_full} shows five system trajectories closely approximating the {\color{blue} set $\mathcal{Y}$} (shown in black). These trajectories result from the full-order model being subjected to slowly time-varying loads and sufficiently slow switching between the loads and the DEA parameters. 

Lastly, Figure \ref{fig:actual_full_time} shows the trajectories versus time of the system subject to the disturbances we mentioned in the previous paragraph. In the top left plot, we see the evolution of the frequencies of the CGs and the DEAs. The frequency deviations remain close to zero with occasional perturbations corresponding to the switching signal (shown in the bottom plot). On the top right, we see the evolution of the ten mechanical powers. Similar in behavior to the frequency deviations, the mechanical powers oscillate among their equilibrium points corresponding to the switching signal. In the top left plot, the bottom inset shows the transient behavior of the frequency deviations shortly after initialization before the behavior starts approximating the $\Omega$-limit set. The inset on the top shows the behavior following a switching event. In this case, we choose the switching event occurring at $t= 500$ seconds (mode $1$ to mode $3$). The deviations converge to zero after a brief transient. 

\vspace{-1mm}

\begin{figure}[t!]
    \centering
    \includegraphics[width=0.9\linewidth]{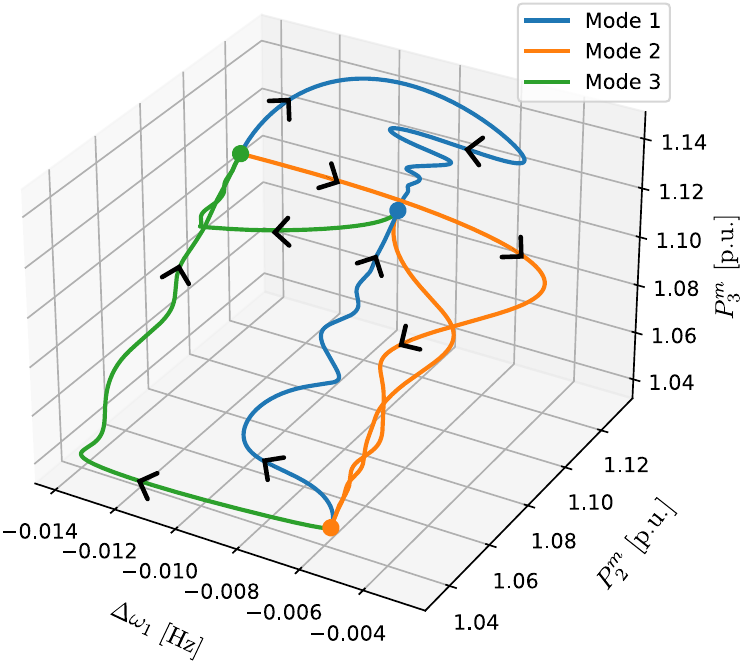}
    \caption{Projection onto the $(\Delta \omega_1, P^m_2, P^m_3)$ axis of the set $\mathcal{Y}$ for the full-order model of the 39-bus system.}
    \label{fig:nominal_39_full}
\end{figure}
\begin{figure}[t!]
    \centering
    \includegraphics[width=0.88\linewidth]{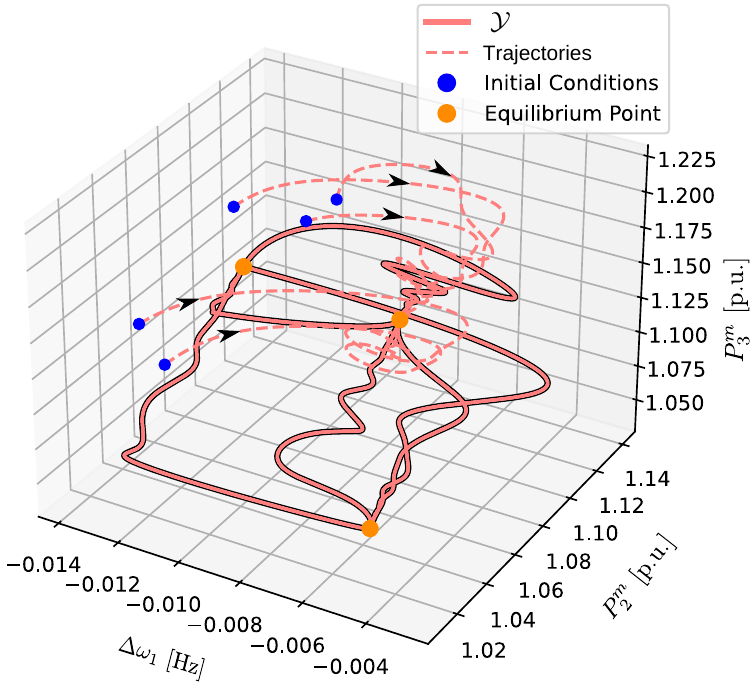}
    \caption{ Illustration of Theorem  1 for the full model of the IEEE 39-bus system.}
    \label{fig:actual_39_full}
\end{figure}
\begin{figure}[t!]
    \centering
    \includegraphics[width=\linewidth]{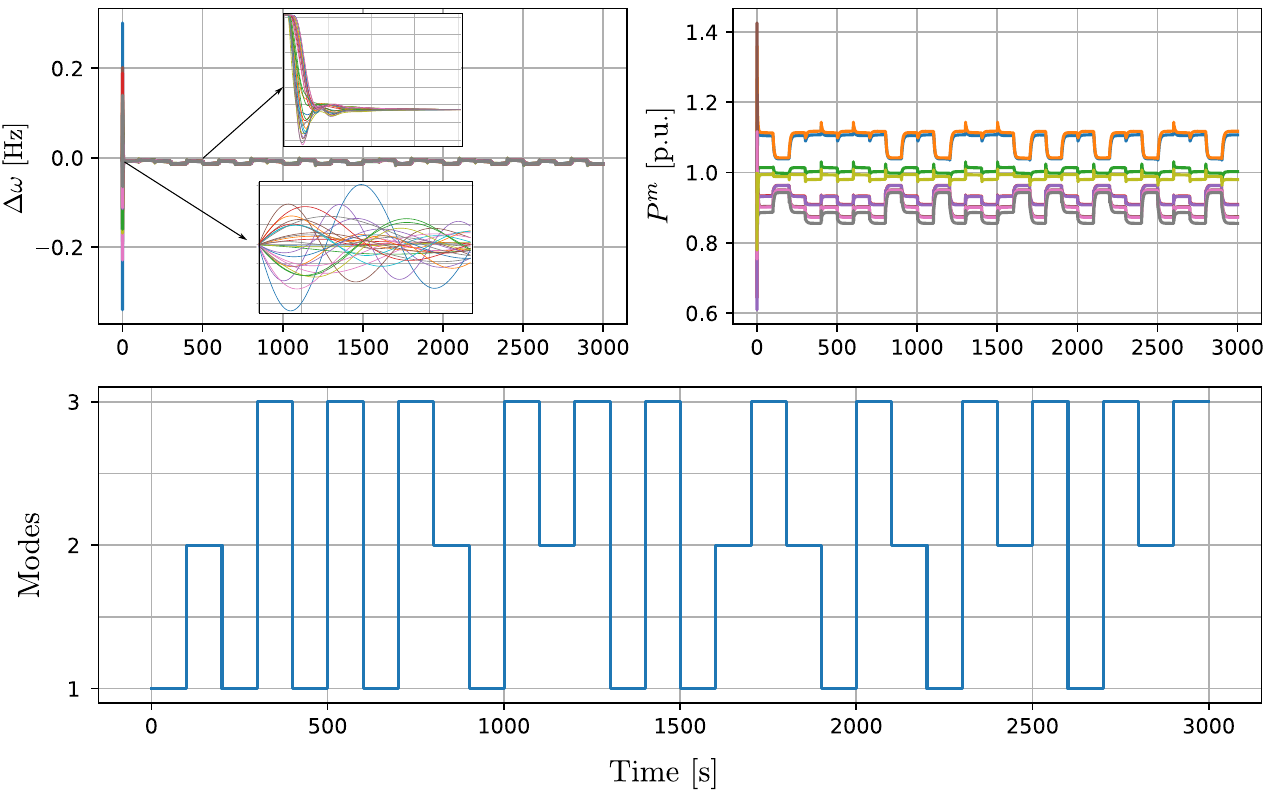}
    \caption{Trajectories versus time of the full-order 39-bus system.}
    \label{fig:actual_full_time}
\end{figure}

\phantom{.}


\phantom{.}

\vspace{-1cm}
\section{Conclusion}
\label{ref:conc}
We modeled the failure states and the presence of intermittent loads for DEAs in a power transmission system via a reduced order model with multiple equilibria and analyzed its asymptotic behavior using the notion of $\Omega$-limit sets. To ease computation for stability analysis in practice, an explicit formula for computing the set in terms of the flows of the subsystems was proved using inner and outer approximations of the reachable sets of a suitably defined nominal system. To deal with time-varying demand profiles, a practical stability result for the case of slowly-varying external disturbances was also proved. Numerical experiments on the IEEE 39-bus system showed that the trajectories of the power transmission system remain close to the $\Omega$-limit set, validating our theoretical results. 

\vspace{2cm}
{\color{black} Future work will extend our analysis to power systems in which some distributed generation units exhibit zero damping $(D_d)$ and inertia $(M_d)$ parameters, leading to algebraic constraints. Such systems can be represented by switching differential-algebraic equations (DAE) with multiple equilibria, the stability of which neccesitates new theoretical developments. We also plan to generalize our framework to multi-area systems by incorporating a network susceptance matrix. Finally, we aim to model demand variations as stochastic processes and reproduce our stability analysis within the framework of stochastic hybrid dynamical systems.}

\section{Appendix}

We present some auxiliary definitions and notions, borrowed from \cite{goebel2012hybrid} and \cite{cai2013robust}, and used throughout this paper. 

\vspace{0.1cm}
A set $E\subset\mathbb{R}_{\geq0}\times\mathbb{Z}_{\geq0}$ is called a \textsl{compact} hybrid time domain if $E=\cup_{j=0}^{J-1}([t_j,t_{j+1}],j)$ for some finite sequence of times $0=t_0\leq t_1\ldots\leq t_{J}$. A set $E \subset \mathbb{R}_{\geq0} \times \mathbb{Z}_{\geq0}$ is a hybrid time domain if it is the union of a non-decreasing sequence of compact hybrid time domains, namely, E is the union of compact hybrid time domains $E_j$ with the property that $E_0 \subset E_1 \subset E_2 \subset \ldots \subset E_j \ldots$, etc. For a given hybrid time domain $E$, and a given $T>0$, let $J^*(T) \coloneqq \inf \{j \, : \, (T,j) \in E\}$. That is, the function $J^*(\cdot)$ returns the first valid jump for which $(T, J^*(T))$ is a point in $E$. Given $(t,j)$ and $(t',j') \in E$, we say $(t,j) \preceq (t',j')$ if $t<t'$ or $t = t'$ and $j \leq j'$. The ordering $\succeq$ is defined analogously.

\vspace{0.1cm}
\begin{defn}
    A function $u:E\to \mathbb{R}^m$ is a \textit{hybrid signal} if $E$ is a hybrid time domain. A hybrid signal is a \textit{hybrid input} if  for each $j\in\mathbb{Z}_{\geq 0}$, the function $t\rightarrow u(t,j)$ is measurable and locally essentially bounded on the interval $I^j = \{t \, : \, (t,j) \in E\}$. A hybrid signal $\phi : E \to \mathbb{R}^n$ is a \textit{hybrid arc} if $E$ is a hybrid time domain and if for each $j \in \mathbb{Z}_{\geq 0}$, the function $t \to \phi(t,j)$ is locally absolutely continuous on the interval $I^j = \{t \, : \, (t,j) \in E\}$. The hybrid time domain of a hybrid signal $\phi$ is denoted as $\text{dom } \phi$. 
\end{defn}

\vspace{0.1cm}
\begin{defn}
    A hybrid arc $\phi$ is a solution to the hybrid system $\mathcal{H}$ \eqref{eqn:hy_sys_def} if $\phi(0,0) \in \Bar{C} \cup D$, and 
    \begin{enumerate}[(S1)]
        \item $\phi(0,0) \in \Bar{C} \cup D$, and $\mathrm{dom}(\phi)=\mathrm{dom}(u)$;
        \item for all $j \in \mathbb{Z}_{\geq 0}$ such that $I^j \coloneqq \{t \, : \, (t,j) \in \text{dom } \phi\}$ has nonempty interior, $\phi(t,j) \in C, ~\text{for all } t \in \text{int } I^j$, and $\dot{\phi}(t,j) \in F(\phi(t,j)), ~\text{for almost all } t \in I^j$;
        \item for all $(t,j) \in \text{dom } \phi$ such that $(t,j+1) \in \text{dom } \phi$, $\phi(t,j) \in D,$ and $\phi(t,j+1) \in G(\phi(t,j)).$
    \end{enumerate}
\end{defn}
\vspace{0.1cm}
\begin{defn}
    A hybrid arc $\phi$ and a hybrid input $u$ form a solution to the hybrid system with input $\mathcal{H}_u$ \eqref{eqn:hy_sys_def_input} if 
    \begin{enumerate}[(S1)]
        \item $\phi(0,0) \in \Bar{C} \cup D$, $\mathrm{dom}(\phi)=\mathrm{dom}(u)$, and $u(t,j)\in\mathcal{U},$ for all $(t,j) \in \text{dom }(u)$;
        \item for all $j \in \mathbb{Z}_{\geq 0}$ such that $I^j \coloneqq \{t \, : \, (t,j) \in \text{dom } \phi\}$ has nonempty interior, $\phi(t,j) \in C, ~\text{for all } t \in \text{int } I^j$, and $\dot{\phi}(t,j) \in F(\phi(t,j),u(t,j)), ~\text{for almost all } t \in I^j$;
        \item for all $(t,j) \in \text{dom } \phi$ such that $(t,j+1) \in \text{dom } \phi$, $\phi(t,j) \in D,$ and $\phi(t,j+1) \in G(\phi(t,j),u(t,j)).$
    \end{enumerate}
\end{defn}

\printcredits

\bibliographystyle{cas-model2-names}

\bibliography{refs}

@article{ochoa2023control,
  title={Control systems for low-inertia power grids: A survey on virtual power plants},
  author={Ochoa, Daniel E and Galarza-Jimenez, Felipe and Wilches-Bernal, Felipe and Schoenwald, David A and Poveda, Jorge I},
  journal={IEEE Access},
  volume={11},
  pages={20560--20581},
  year={2023},
  publisher={IEEE}
}

@article{zografopoulos2021detection,
  title={Detection of malicious attacks in autonomous cyber-physical inverter-based microgrids},
  author={Zografopoulos, Ioannis and Konstantinou, Charalambos},
  journal={IEEE Transactions on Industrial Informatics},
  volume={18},
  number={9},
  pages={5815--5826},
  year={2021},
  publisher={IEEE}
}

@inproceedings{barua2020hall,
  title={Hall {S}poofing: A {N}on-{I}nvasive {D}o{S} Attack on {G}rid-{T}ied Solar Inverter},
  author={Barua, Anomadarshi and Al Faruque, Mohammad Abdullah},
  booktitle={29th USENIX Security Symposium (USENIX Security 20)},
  pages={1273--1290},
  year={2020}
}

@inproceedings{kang2015investigating,
  title={Investigating cyber-physical attacks against {IEC} 61850 photovoltaic inverter installations},
  author={Kang, BooJoong and Maynard, Peter and McLaughlin, Kieran and Sezer, Sakir and Andr{\'e}n, Filip and Seitl, Christian and Kupzog, Friederich and Strasser, Thomas},
  booktitle={IEEE Conference on Emerging Technologies \& Factory Automation (ETFA)},
  pages={1--8},
  year={2015},
}

@inproceedings{wheeler2018power,
  title={A power quality and load analysis of a cryptocurrency mine},
  author={Wheeler, Keaton A and Bowers, Anthony W and Wong, Charlie H and Palmer, Jonathan Y and Wang, Xue},
  booktitle={2018 IEEE Elec. Power and Energy Conf. (EPEC)},
  pages={1--6},
  year={2018},
  organization={IEEE}
}

@article{goebel2012hybrid,
  title={Hybrid dynamical systems: {M}odeling, {S}tability, and {R}obustness},
  author={Goebel, Rafal and Sanfelice, Ricardo G and Teel, Andrew R},
  journal={Princeton, NJ, USA},
  year={2012}
}

@book{kundur2007power,
	title = {Power {System} {Stability} and {Control}},
	publisher = {
    McGraw-Hill},
	author = {Kundur, Prabha},
    year = {1994}
}

@article{dorfler2023control,
  title={Control of low-inertia power systems},
  author={D{\"o}rfler, Florian and Gro{\ss}, Dominic},
  journal={Annual Review of Control, Robotics, and Autonomous Systems},
  volume={6},
  pages={415--445},
  year={2023},
  publisher={Annual Reviews}
}

@article{mavalizadeh2023improving,
  title={Improving frequency response with synthetic damping available from fleets of distributed energy resources},
  author={Mavalizadeh, Hani and Espinosa, Luis A Duffaut and Almassalkhi, Mads R},
  journal={IEEE Trans. on Power Syst.},
  year={2023},
  publisher={IEEE}
}

@article{cheng2020smart,
  title={Smart frequency control in low inertia energy systems based on frequency response techniques: A review},
  author={Cheng, Yi and Azizipanah-Abarghooee, Rasoul and Azizi, Sadegh and Ding, Lei and Terzija, Vladimir},
  journal={Applied Energy},
  volume={279},
  pages={115798},
  year={2020},
  publisher={Elsevier}
}

@article{anderson1990low,
  title={A low-order system frequency response model},
  author={Anderson, Philip M and Mirheydar, Mahmood},
  journal={IEEE Trans. on Power Syst.},
  volume={5},
  number={3},
  pages={720--729},
  year={1990},
  publisher={IEEE}
}

@book{wood2013power,
  title={Power generation, operation, and control},
  author={Wood, Allen J and Wollenberg, Bruce F and Shebl{\'e}, Gerald B},
  year={2013},
  publisher={John Wiley \& Sons}
}

@article{guggilam2018optimizing,
  title={Optimizing {DER} participation in inertial and primary-frequency response},
  author={Guggilam, Swaroop S and Zhao, Changhong and Dall’Anese, Emiliano and Chen, Yu Christine and Dhople, Sairaj V},
  journal={IEEE Trans. on Power Syst.},
  volume={33},
  number={5},
  pages={5194--5205},
  year={2018},
  publisher={IEEE}
}

@article{colombino2019online,
  title={Online optimization as a feedback controller: Stability and tracking},
  author={Colombino, Marcello and Dall’Anese, Emiliano and Bernstein, Andrey},
  journal={IEEE Transactions on Control of Network Systems},
  volume={7},
  number={1},
  pages={422--432},
  year={2019},
  publisher={IEEE}
}

@article{baradaran2020omega,
  title={Omega-limit sets and robust stability for switched systems with distinct equilibria},
  author={Baradaran, Matina and Teel, Andrew R},
  journal={IFAC-PapersOnLine},
  volume={53},
  number={2},
  pages={2039--2044},
  year={2020},
  publisher={Elsevier}
}

@article{alpcan2010stability,
  title={A stability result for switched systems with multiple equilibria},
  author={Alpcan, Tansu and Basar, Tamer},
  journal={Dynamics of Continuous, Discrete and Impulsive Systems Series A: Mathematical Analysis},
  volume={17},
  number={4},
  pages={949--958},
  year={2010},
  publisher={Watam Press}
}

@article{yin2023stability,
  title={Stability of switched systems with multiple equilibria: A mixed stable--unstable subsystem case},
  author={Yin, Hao and Jayawardhana, Bayu and Trenn, Stephan},
  journal={Systems \& Control Letters},
  volume={180},
  pages={105622},
  year={2023},
  publisher={Elsevier}
}

@article{veer2019switched,
  title={Switched systems with multiple equilibria under disturbances: Boundedness and practical stability},
  author={Veer, Sushant and Poulakakis, Ioannis},
  journal={IEEE Transactions on Automatic Control},
  volume={65},
  number={6},
  pages={2371--2386},
  year={2019},
  publisher={IEEE}
}

@article{horowitz2019overview,
  title={An overview of distributed energy resource ({DER}) interconnection: Current practices and emerging solutions},
  author={Horowitz, Kelsey A and Peterson, Zachary and Coddington, Michael H and Ding, Fei and Sigrin, Benjamin O and Saleem, Danish and Baldwin, Sara E and Lydic, Brian and Stanfield, Sky C and Enbar, Nadav and others},
  year={2019},
  journal={},
  publisher={National Renewable Energy Lab.(NREL), Golden, CO (United States)}
}

@article{dibaji2019systems,
  title={A systems and control perspective of CPS security},
  author={Dibaji, Seyed Mehran and Pirani, Mohammad and Flamholz, David Bezalel and Annaswamy, Anuradha M and Johansson, Karl Henrik and Chakrabortty, Aranya},
  journal={Annual reviews in control},
  volume={47},
  pages={394--411},
  year={2019},
  publisher={Elsevier}
}

@article{kavvathas2024resilient,
  title={Resilient Distributed Control for Power Systems With Multiple Synchronous Generators and {DER}s Satisfying Capability Curve Requirements},
  author={Kavvathas, Theodoros and Konstantopoulos, George C},
  journal={IEEE Trans on Power Syst},
  year={2024},
  publisher={IEEE}
}

@article{tang2023common,
  title={Common {L}yapunov function based stability analysis of VSC with limits of phase locked loop},
  author={Tang, Yingyi and Li, Yujun},
  journal={IEEE Trans on Power Syst.},
  volume={38},
  number={2},
  pages={1759--1762},
  year={2023},
  publisher={IEEE}
}

@article{feng2025online,
  title={Online Event-Triggered Switching for Frequency Control in Power Grids With Variable Inertia},
  author={Feng, Jie and Cui, Wenqi and Cort{\'e}s, Jorge and Shi, Yuanyuan},
  journal={IEEE Trans on Power Syst.},
  year={2025},
  publisher={IEEE}
}

@article{chen2023distributed,
  title={Distributed dual objective control of an energy storage system under jointly connected switching networks},
  author={Chen, Xiaofeng and Cai, He and Su, Youfeng},
  journal={Automatica},
  volume={152},
  pages={110979},
  year={2023},
  publisher={Elsevier}
}

@article{peng2022distributed,
  title={Distributed periodic event-triggered optimal control of {DC} microgrids based on virtual incremental cost},
  author={Peng, Jiangkai and Fan, Bo and Tu, Zhenghong and Zhang, Wei and Liu, Wenxin},
  journal={IEEE/CAA Journal of Automatica Sinica},
  volume={9},
  number={4},
  pages={624--634},
  year={2022},
  publisher={IEEE}
}

@article{standard,
    author = {{CENELEC}},
    title = {Requirements for Generating Plants to be Connected in Parallel With {D}istribution Networks- {P}art 2 : Connection
            to a {MV} {D}istribution Network - {G}enerating Plant up to and Including {T}ype {B},},
    journal ={CENELEC, Brussels, Belgium,} ,
    year = {2019}
}

@article{ghiasi2019analytical,
  title={An analytical methodology for reliability assessment and failure analysis in distributed power system},
  author={Ghiasi, Mohammad and Ghadimi, Noradin and Ahmadinia, Esmaeil},
  journal={SN Applied Sciences},
  volume={1},
  number={1},
  pages={44},
  year={2019},
  publisher={Springer}
}

@inproceedings{jung2024inferences,
  title={Inferences from Numerical Model Reduction \& Aggregation of Power-system Dynamics Featuring Primary \& Secondary Control},
  author={Jung, HyeonjungTari and Helman, Samuel and Singh, Manish K and Dhople, Sairaj},
  booktitle={2024 IEEE Power \& Energy Society General Meeting (PESGM)},
  pages={1--5},
  year={2024},
  organization={IEEE}
}

@article{wang2024distributed,
  title={Distributed Coordination of Grid-Forming and Grid-Following Inverter-Based Resources for Optimal Frequency Control in Power Systems},
  author={Wang, Xiaoyang and Chen, Xin},
  journal={arXiv preprint arXiv:2411.12682},
  year={2024}
}

@article{ochoa2024prescribed,
  title={Prescribed-time stability in switching systems with resets: A hybrid dynamical systems approach},
  author={Ochoa, Daniel E and Espitia, Nicolas and Poveda, Jorge I},
  journal={Systems \& Control Letters},
  volume={193},
  pages={105910},
  year={2024},
  publisher={Elsevier}
}

@article{cai2009characterizations,
  title={Characterizations of input-to-state stability for hybrid systems},
  author={Cai, Chaohong and Teel, Andrew R},
  journal={Systems \& Control Letters},
  volume={58},
  number={1},
  pages={47--53},
  year={2009},
  publisher={Elsevier}
}

@article{zhou2013stability,
  title={Stability analysis of switched linear singular systems},
  author={Zhou, Lei and Ho, Daniel WC and Zhai, Guisheng},
  journal={Automatica},
  volume={49},
  number={5},
  pages={1481--1487},
  year={2013},
  publisher={Elsevier}
}

@article{cai2013robust,
  title={Robust input-to-state stability for hybrid systems},
  author={Cai, Chaohong and Teel, Andrew R},
  journal={SIAM Journal on Control and Optimization},
  volume={51},
  number={2},
  pages={1651--1678},
  year={2013},
  publisher={SIAM}
}

@article{liberzon2012switched,
  title={Switched nonlinear differential algebraic equations: Solution theory, {L}yapunov functions, and stability},
  author={Liberzon, Daniel and Trenn, Stephan},
  journal={Automatica},
  volume={48},
  number={5},
  pages={954--963},
  year={2012},
  publisher={Elsevier}
}

@article{debeljkovic2011stability,
  title={Stability of linear discrete descriptor systems in the sense of {L}yapunov},
  author={Debeljkovic, D and Buzurovic, I and Simeunovic, G},
  journal={International journal of information and systems sciences},
  volume={7},
  number={4},
  pages={303--322},
  year={2011}
}

@standard{IEEE-Std-1547-2018,
  title        = {{IEEE Standard for Interconnection and Interoperability of Distributed Energy Resources with Associated Electric Power Systems Interfaces}},
  author       = {IEEE standards association},
  number       = {IEEE Std 1547-2018},
  year         = {2018},
  doi          = {10.1109/IEEESTD.2018.8332112},
  url          = {https://standards.ieee.org/standard/1547-2018.html},
  note         = {Interconnection and performance requirements for DERs}
}

@book{kokotovic1999singular,
  title={Singular perturbation methods in control: analysis and design},
  author={Kokotovi{\'c}, Petar and Khalil, Hassan K and O'{R}eilly, John},
  year={1999},
  publisher={SIAM}
}

@book{chow1982time,
  title={Time-scale modeling of dynamic networks with applications to power systems},
  author={Chow, Joe H},
  year={1982},
  publisher={Springer}
}

@ARTICLE{Yuantime,
  author={Yuan, Hao and Yuan, Xiaoming and Hu, Jiabing},
  journal={IEEE Transactions on Power Systems}, 
  title={Modeling of Grid-Connected VSCs for Power System Small-Signal Stability Analysis in DC-Link Voltage Control Timescale}, 
  year={2017},
  volume={32},
  number={5},
  pages={3981-3991},
  doi={10.1109/TPWRS.2017.2653939}}

@article{ajala2022model,
  title={Model reduction and dynamic aggregation of grid-forming inverter networks},
  author={Ajala, Olaoluwapo and Baeckeland, Nathan and Johnson, Brian and Dhople, Sairaj and Dom{\'\i}nguez-Garc{\'\i}a, Alejandro},
  journal={IEEE Transactions on Power Systems},
  volume={38},
  number={6},
  pages={5475--5490},
  year={2022},
  publisher={IEEE}
}

@article{he2020transient,
  title={Transient stability of power systems integrated with inverter-based generation},
  author={He, Xiuqiang and Geng, Hua},
  journal={IEEE Transactions on power systems},
  volume={36},
  number={1},
  pages={553--556},
  year={2020},
  publisher={IEEE}
}

@article{cotilla2013multi,
  title={Multi-attribute partitioning of power networks based on electrical distance},
  author={Cotilla-Sanchez, Eduardo and Hines, Paul DH and Barrows, Clayton and Blumsack, Seth and Patel, Mahendra},
  journal={IEEE Transactions on Power Systems},
  volume={28},
  number={4},
  pages={4979--4987},
  year={2013},
  publisher={IEEE}
}

@book{chow2013power,
  title={Power system coherency and model reduction},
  author={Chow, Joe H},
  volume={84},
  year={2013},
  publisher={Springer}
}

@article{odunlami2025hybrid,
  title={Hybrid dynamical systems modeling of power systems},
  author={Odunlami, BG and Netto, M and Susuki, Y},
  journal={arXiv preprint arXiv:2509.02822},
  year={2025}
}

@article{susuki2011hybrid,
  title={A hybrid system approach to the analysis and design of power grid dynamic performance},
  author={Susuki, Yoshihiko and Koo, T John and Ebina, Hiroaki and Yamazaki, Takuya and Ochi, Takashi and Uemura, Takuji and Hikihara, Takashi},
  journal={Proceedings of the IEEE},
  volume={100},
  number={1},
  pages={225--239},
  year={2011},
  publisher={IEEE}
}

@article{samanta2025stability,
  title={Stability analysis of power systems with high penetration of state-of-the-art inverter technologies},
  author={Samanta, Sayan and Yang, Bowen and Seo, Gab-Su},
  journal={Energies},
  volume={18},
  number={14},
  pages={3645},
  year={2025},
  publisher={MDPI}
}

@article{bakhshizadeh2023improved,
  title={Improved reduced-order model for {PLL} instability investigations},
  author={Bakhshizadeh, Mohammad Kazem and Ghosh, Sujay and Kocewiak, {\L}ukasz and Yang, Guangya},
  journal={IEEE Access},
  volume={11},
  pages={72400--72408},
  year={2023},
  publisher={IEEE}
}

@article{schiffer2016survey,
  title={A survey on modeling of microgrids—From fundamental physics to phasors and voltage sources},
  author={Schiffer, Johannes and Zonetti, Daniele and Ortega, Romeo and Stankovi{\'c}, Aleksandar M and Sezi, Tevfik and Raisch, J{\"o}rg},
  journal={Automatica},
  volume={74},
  pages={135--150},
  year={2016},
  publisher={Elsevier}
}

@book{sanfelice2021hybrid,
  title={Hybrid feedback control},
  author={Sanfelice, Ricardo G},
  year={2021},
  publisher={Princeton University Press}
}

\end{document}